\documentclass[a4paper,12pt]{article}
\usepackage{dsfont}
\usepackage{amssymb}
\usepackage{latexsym}
\usepackage{amsmath}
\usepackage{color}
\usepackage{comment}
\usepackage{amsthm}
\usepackage{graphicx}
\usepackage{layout} 
\usepackage{amsfonts}
\usepackage{amsmath}
\usepackage{ulem}
\usepackage{enumerate}
\usepackage{bm}
\usepackage{bbm} 
\usepackage[bbgreekl]{mathbbol} 
\usepackage{authblk}
\usepackage{setspace} 
\usepackage{here}
\newif\ifrs
\rstrue
\ifrs \usepackage{mathrsfs} \fi  
\newif\ifcol
\coltrue 

\newtheorem{theorem}{Theorem}[section]

\newtheorem{lemma}[theorem]{Lemma}

\newtheorem{remark}[theorem]{Remark}
\newtheorem{example}[theorem]{Example}
\numberwithin{equation}{section}
\newtheorem{theorem*}{Theorem}
\newtheorem{ass*}[theorem*]{Assumption}
\newtheorem{note*}[theorem*]{Note}
\newtheorem{lemma*}[theorem*]{Lemma}
\newtheorem{definition*}[theorem*]{Definition}
\newtheorem{proposition*}[theorem*]{Proposition}
\newtheorem{corollary*}[theorem*]{Corollary}
\newtheorem{remark*}[theorem*]{Remark}
\newtheorem{example*}[theorem*]{Example}
\numberwithin{equation}{section}
\newcommand\ttU{{\tt U}}
\newcommand\bbd{{\mathbbm d}}
\newcommand\consta{z}
\newcommand\constb{C_1}
\newcommand\constd{C_{23}}
\newcommand\constg{C_2}
\newcommand\consth{{C_3}}
\newcommand\consti{{C_4}}
\newcommand\constj{{C_6}}
\newcommand\constk{{C_5}}
\newcommand\constkp{{C_5'}}
\newcommand\constl{{C_7}}
\newcommand\constn{{C_9}}
\newcommand\consto{{C_8}}
\newcommand\constp{{C_{19}}}
\newcommand\constq{{C_{0}}}
\newcommand\constr{{C_{12}}}
\newcommand\consts{{C_{16}}}
\newcommand\constt{{C_{18}}}
\newcommand\constu{{C_{17}}}
\newcommand\constv{{C_{10}}}
\newcommand\constw{{C_{11}}}
\newcommand\constx{{C_{20}}}
\newcommand\constab{{C_{13}}}
\newcommand\constac{{C_{15}}}
\newcommand\constad{{C_{14}}}
\newcommand\constzz{{C_{*}}}
\newcommand\ttT{{\tt T}}
\newcommand\tth{{\tt h}}
\newif\ifcol
\colfalse
\ifcol
\newcommand{\colorr}{\color[rgb]{0.8,0,0}}

\newcommand{\colorb}{\color[rgb]{0,0,0.8}}

\newcommand{\colorn}{\color[rgb]{1,1,1}}

\newcommand{\colory}{\color{yellow}}

\else
\newcommand{\colorb}{\color{black}}

\newcommand{\colorr}{\color{black}}

\newcommand{\colorn}{\color{black}}
\newcommand{\colory}{\color{black}}

\fi
\newif\ifcol
\colfalse
\ifcol
\newcommand{\cred}{\color[rgb]{0.8,0,0}}

\else
\newcommand{\cred}{\color{black}}
\fi
\newif\ifcol
\colfalse
\ifcol

\newcommand{\tred}{\color[rgb]{0.8,0,0}}

\else
\newcommand{\tred}{\color{black}}
\fi
\newif\ifcol
\colfalse
\ifcol

\else
\fi
\newif\ifcol
\colfalse
\ifcol

\else
\fi
\newif\ifcol
\colfalse
\ifcol

\newcommand{\ared}{\color[rgb]{0.8,0,0}}

\newcommand{\ablue}{\color[rgb]{0,0,0.8}}

\else
\newcommand{\ared}{\color{black}}
\newcommand{\ablue}{\color{black}}
\fi
\newif\ifcol
\colfalse
\ifcol

\newcommand{\bred}{\color[rgb]{0.8,0,0}}

\newcommand{\bblue}{\color[rgb]{0,0,0.8}}

\else
\newcommand{\bred}{\color{black}}
\newcommand{\bblue}{\color{black}}
\fi
\excludecomment{en-text}
\includecomment{jp-text}
\includecomment{comment}
\def\koko{{\coloroy{koko}}}
\def\bd{\begin{description}}
\def\ed{\end{description}}

\def\D2{\bbD_{2,\infty-}}

\def\D{{\bf D}}

\def\calb{{\cal B}}

\def\cald{{\cal D}}

\def\calf{{\cal F}}
\def\calg{{\cal G}}

\def\calk{{\cal K}}

\def\caln{{\cal N}}

\def\calr{{\cal R}}

\def\calv{{\cal V}}

\def\calx{{\cal X}}

\def\yeq{\>=\>}
\def\yleq{\>\leq\>}
\def\ygeq{\>\geq\>}

\def\sfd{{\sf d}}
\def\sfp{{\sf p}}

\def\simleq{\ \raisebox{-.7ex}{$\stackrel{{\textstyle <}}{\sim}$}\ }

\def\ep{\epsilon}
\def\half{\frac{1}{2}}

\def\down{\downarrow}

\def\halflineskip{\vspace*{3mm}}
\def\nn{\nonumber}
\def\be{\begin{equation}}
\def\ee{\end{equation}}
\def\bea{\begin{eqnarray}}
\def\eea{\end{eqnarray}}
\def\beas{\begin{eqnarray*}}
\def\eeas{\end{eqnarray*}}
\def\bi{\begin{itemize}}
\def\ei{\end{itemize}}
\def\im{\item}
\def\bd{\begin{description}}
\def\ed{\end{description}}
\newcommand{\bbA}{{\mathbb A}}
\newcommand{\bbB}{{\mathbb B}}

\newcommand{\bbD}{{\mathbb D}}

\newcommand{\bbI}{{\mathbb I}}

\newcommand{\bbK}{{\mathbb K}}
\newcommand{\bbL}{{\mathbb L}}
\newcommand{\bbM}{{\mathbb M}}
\newcommand{\bbN}{{\mathbb N}}

\newcommand{\bbR}{{\mathbb R}}

\newcommand{\bbT}{{\mathbb T}}

\newcommand{\bbX}{{\mathbb X}}

\newcommand{\bbZ}{{\mathbb Z}}

\newcommand{\ol}{\overline}

\renewcommand{\colory}{\color{yellow}}
\renewcommand{\koko}{{\colory koko}}
\newcommand{\wh}{\widehat}
\newcommand{\wt}{\widetilde}

\newcommand{\gray}{\color[rgb]{0.5,0.5,0.5}}

\begin{document}

\title{
Deep learning of point processes for modeling high-frequency data
\footnote{
This work was in part supported by 
Japan Science and Technology Agency CREST JPMJCR2115;  
Japan Society for the Promotion of Science Grants-in-Aid for Scientific Research 
No. 23H03354 (Scientific Research);  
Forefront Physics and Mathematics Program to Drive Transformation (FoPM), a World-leading Innovative Graduate Study (WINGS) Program, the University of Tokyo; 
and by a Cooperative Research Program of the Institute of Statistical Mathematics. 
}
}
\author[1,3]{Yoshihiro Gyotoku\thanks{gyotoku@ms.u-tokyo.ac.jp}}
\author[2,3]{Ioane Muni Toke\thanks{ioane.muni-toke@centralesupelec.fr}}
\author[1,3]{Nakahiro Yoshida\thanks{nakahiro@ms.u-tokyo.ac.jp}}
\affil[1]{University of Tokyo, Graduate School of Mathematical Sciences
\footnote{Graduate School of Mathematical Sciences, University of Tokyo: 3-8-1 Komaba, Meguro-ku, Tokyo 153-8914, Japan. }
        }
\affil[2]{{\ablue 
Universit\'e Paris-Saclay, CentraleSup\'elec, Math\'ematiques et Informatique pour la Complexit\'e et les Syst\`emes
}
\footnote{{\ablue Universit\'e Paris-Saclay, CentraleSupelec, 3 rue Joliot Curie, 91190 Gif-sur-Yvette, France}}
}
\affil[3]{Japan Science and Technology Agency CREST
       }
\maketitle
\ \\
{\it Summary} 
{\ared
We investigate applications of deep neural networks to a point process having an intensity with mixing covariates processes as input. Our generic model includes Cox-type models and marked point processes as well as multivariate point processes. 
An oracle inequality and a rate of convergence are derived for the prediction error. {\ablue A} simulation study shows that the marked point process can be superior to the simple multivariate model in prediction. We apply the marked ratio model to {\ablue real limit order book data}. 
}



\section{Point process with covariates and deep learning}
Given a stochastic basis $\calb=(\Omega,\calf,{\bf F},P)$, 
we consider a $\sfd_N$-dimensional counting process $N=(N^i)_{i\in\bbI}$ 
with a $\sfd_N$-dimensional intensity process {\ared $\lambda_t$} with respect to ${\bf F}$, where $\bbI$ is a finite index set with $\#\bbI=\sfd_N$. 
It is assumed that 
${\bf F}=(\calf_t)_{t\in\bbR_+}$ is a right-continuous filtration and the stochastic basis $\calb$ satisfies the usual condition. 
The intensity process {\ared $\lambda_t$} is supposed to admit a representation 
{\ared $\lambda_t=\lambda(X_t)$} with 
a bounded measurable mapping ${\ared \lambda=((\lambda^i))_{i\in\bbI}}:\calx\to\bbR^{\sfd_N}$ and 
a $\sfd_X$-dimensional ${\bf F}$-predictable covariate process $X=(X_t)_{t\in\bbR_+}$ 
taking values in $\calx\in\bbB[\bbR^{\sfd_X}]$, the Borel $\sigma$-field. 
{\ared Suppose that the true mechanism that generates the data $N$ is denoted by a mapping $\lambda^*$ among possible mappings $\lambda$. }%
{\bred Moreover, we suppose that the components $N^i$ have no common jumps. }%

For a positive number ${\tt h}$, let 
$I_j=((j-1){\tt h},j{\tt h})$ ($j\in\bbN=\{1,2,...\}$) and $\bbX_j=(X_t,N_t-N_s)_{t,s\in I_j}$. 
Suppose that $(X,N)$ is periodically stationary, that is, $(\bbX_j)_{j\in\bbN}$ is stationary. 
For example, the periodical stationarity models the stochastic evolution of a limit order book which has intraday non-stationarity but has a long-term stationarity. 
We will consider the periodically stationary case in this paper, while the (exactly) stationary case can be dealt with, as a matter of fact, more simply.

{\ared The model of the intensity process is expressed by $\lambda:\calx\to\bbR^{\sfd_N}$, as already mentioned. }%
Let $T\in\bbT={\tt h}\bbN$. 
Consider certain functions $a=(a_i)_{i\in\bbI}$ and $b$ of $\lambda$, more general than $\lambda$ itself; 
see the examples below. 
Let $(a^*,b^*)$ denote $(a,b)$ for $\lambda^*$.
We estimate the mapping $(a^*,b^*)$ from the data $(X_t,N_t)_{t\in[0,T]}$ 
with a family ${\mathfrak F}_T$ of candidates mappings $(a,b)$.  
It is not assumed that $(a^*,b^*)$ belongs to the family ${\mathfrak F}_T$. 
Three examples of setting $(a,b)$ will be provided later. 
A family ${\mathfrak F}_T$ we are interested in in this article is a deep neural network the size of which is increasing to infinity as $T\to\infty$. 
However, the result we will obtain is more general and not confined to {\ablue the case of deep learning (DL)}. 
The aim of this paper is to derive a bound for the prediction error by the machine ${\mathfrak F}_T$ 
applied to point processes. 

{\ared 

Modeling limit order book (LOB) with point processes has been a big trend for the past two decades. Early point process models include
Bowsher \cite{bowsher2007modelling}, 
Large \cite{large2007measuring}, 
Bacry et al. \cite{bacry2012non,bacry2013modelling}, 
Muni Toke and Pomponio \cite{muni2011modelling}, 
Lu and Abergel \cite{lu2018high}, just to mention few.
More recent contributions propose intensity models depending on observable LOB covariates:
Muni Toke and Yoshida \cite{muni2017modelling}, 
Rambaldi et al. \cite{rambaldi2017role}, 
Morariu-Patrichi and Pakkanen \cite{morariu2022state}, 
Wu et al. \cite{wu2022queue}, Sfendourakis et al. \cite{sfendourakis2023lob}.
Deep learning architectures have also been proposed for the modeling of limit order book, see e.g. Tsantekidis et al. \cite{tsantekidis2017forecasting}, Sirignano \cite{sirignano2019deep}, Zhang et al. \cite{zhang2019deeplob}, Maglaras and Moallemi \cite{maglaras2022deep} among many others. Many deep learning contributions focus on the prediction of price movements at a given horizon using some specifically designed neural network architecture feeded with limit order book features.
}

%
%

This paper's attempt to incorporate deep learning to point processes is motivated by the authors' studies 
on modeling of the limit order book. 
Muni Toke and Yoshida \cite{muni2019analyzing} took a parametric approach with a Cox-type model 
(the ratio model) 
for relative intensities of order flows in the limit order book. 
The Cox-type model with a nuisance baseline hazard is well suited to cancel 
non-stationary intraday trends in the market data. 
They showed consistency and asymptotic normality of a quasi-likelihood estimator and validated the model selection criteria applied to the point processes, based on the quasi-likelihood analysis (Yoshida \cite{Yoshida2011,yoshida2024simplified}). 
Their scheme is applied to {\ablue real data} from the Paris Stock Exchange 
and achieves accurate prediction of the signs of market orders, as 
the method outperforms the traditional Hawkes model. 
It is suggested that the selection of the covariates is crucial for prediction. 
Succeedingly, Muni Toke and Yoshida \cite{muni2022marked} extended the ratio model to a marked ratio model to express a hierarchical structure in market orders. 
Each market order is categorized by {\ablue bid/ask}, further into aggressive/non-aggressive orders according to the existence of price change. The marked ratio model outperforms other intensity-based methods like Hawkes-based methods in predicting the sign and aggressiveness of market orders on financial markets.
However, the trials of model selection in \cite{muni2019analyzing,muni2022marked} suggest a possibility of taking more covariates in the model; the information criteria seem to prefer relatively large models among a large number of models generated by combinations of our proposal covariates. 
This motivates us to use deep learning to automatically generate more covariates and to enhance the power of expression of the model for more nonlinear dependencies behind the data.

{\ared 
According to a recent big surge of applications of {\ablue deep learning}, theoretical analysis of the prediction error has been a hot topic in the nonparametric statistical approaches to it.
Among these efforts, several survey papers (e.g., \cite{SuhCheng2024survey, Farrell2021deep, DeVore2021neural, Fan2021selective}) provide a comprehensive overview of the state of the art, offering valuable insights into the key open questions and major developments in the field.

More specifically, our work builds on the seminal research presented by Schmidt-Hieber \cite{schmidt2020nonparametric}, which analysed the nonparametric estimation of a specific class of function using fully connected feed forward neural networks with ReLU (Rectified Linear Unit) activation function under independent and identically distributed observations.
Since the publication of Schmidt-Hieber \cite{schmidt2020nonparametric}, several subsequent works \cite{Imaizumi2023supnorm, Kim2024transformers, Kurisu2025adaptive, Oko2023diffusion} have explored its ideas further, extending or applying them to other types of data with various dependency and/or more sophisticated neural network architectures.
}

The organization of this paper is as follows. 
In Section \ref{202412271342}, we formulate the problem more precisely and give a theoretical result on the prediction error, 
whose proof is given in Section \ref{0610080715}. 
The result is not restricted to {\ablue the case of deep learning}. 
Section \ref{202412271344} treats an application of the above result to {\ablue the case of deep learning}. 
The ratio model is investigated in Section \ref{202412271346} in the light of {\ablue deep learning}. 
It will be shown that information on the structure of the model can serve to diminish the error even if it is nonparametric, and that 
it is the case when one uses {\ablue deep learning models}.

\section{Rate of convergence of the error}\label{202412271342}
Let us start the discussion with the loss function defining the prediction error. 
We introduce a contrast function 
\beas 
\Psi_T(a,b) 
&=& 
-\int_0^{T}a(X_t)\cdot dN_t+\int_0^{T}b(X_t)dt\qquad(T\in\bbT). 
\eeas
The discrepancy between $(a,b)$ and $(a^*,b^*)$ is assessed by the function
\beas 
U(x)&=&{\colorr U^{(a,b)}(x)}\yeq
-\lambda^*(x)\cdot\big\{a(x)-a^*(x)\big\}+\big\{b(x)-b^*(x)\big\}. 
\eeas  
Assume that $U(x)\geq0$ for all $x\in\calx$.

\begin{example}\label{202412120701}\rm (Likelihood)  
The {\ablue minus log-likelihood function} $\Psi_T(a,b)$ is realized as 
$\lambda(x)=\big(\lambda^i(x)\big)_{i\in\bbI}$, 
$a_i(x)=\log\lambda^i(x)$ ($i\in\bbI$) and $b(x)=\sum_{i\in\bbI}\lambda^i(x)$. 
Then 
$U(x)\geq0$. 
\end{example}

\begin{example}\label{202504112106}\rm (Ratio model)
{\ared The ratio model of Muni Toke and Yoshida \cite{muni2019analyzing} uses $r^i(x)=\lambda^i(x)/\sum_{i'\in\bbI}\lambda^{i'}(x)$ {\ablue ($i\in\bbI$)}. }%
In this case, 
$a(x)=\big(\log r^i(x)\big)_{i\in\bbI}$ and $b(x)=0$, 
Then 
$U(x)\geq0$. 
As a generalization, 
Muni Toke and Yoshida \cite{muni2022marked} considered a marked ratio model. 
The loss functions for the marked ratio model are exemplified in Section \ref{202412271346}. 
\end{example}

\begin{example}\label{202504112107}\rm
(Mixed loss) 
{\ared
The marks for the $i$-th counting process $N^i$ take values in a finite set $\bbK_i$. 
The process $N^{i,k_i}$ counts the number of the events $(i,k_i)$, and 
the {\ablue intensities are} given by }
\beas 
\lambda^{i,k_i}(X_t,Y_t) &=& \lambda^i(X_t)p_i^{k_i}(Y^i_t),\quad \sum_{k_i\in\bbK_i}p_i^{k_i}\yeq1, 
\eeas
{\ared where $Y^i=(Y^i_t)_{t\in\bbR_+}$ is a covariate process for the mark process associated with $N^i$. }%
Then the likelihood type loss function becomes a mixture of log-likelihoods of a point process and a ratio model:  
\beas &&
-\sum_{i,k_i}\int_0^T\log\big\{\lambda^i(X_t)p_i^{k_i}(Y^i_t)\big\}dN^{i,k_i}_t+\sum_{i,k_i}\int_0^T\lambda^i(X_t)p_i^{k_i}(Y^i_t)dt
\nn\\&=&
-\bigg\{\int_0^T\sum_i\log\lambda^i(X_t)dN^i_t-\int_0^T\sum_i\lambda^i(X_t)dt\bigg\}
-\int_0^T\sum_{i,k_i}\log p_i^{k_i}(Y^i_t)dN^{i,k_i}_t,
\eeas
where $N^i=\sum_{k_i\in\bbK^i}N^{i,k_i}$. 
In this case, 
$a(x,y)=\big((\log\lambda^i(x))_{i\in\bbI},(\log p_i^{k_i}(y))_{i\in\bbI,k_i\in\bbK^i}\big)$ 
and $b(x,y)=\sum_{i\in\bbI}\lambda^i(x)$, 
for the multivariate point process $\big((N^i)_{i\in\bbI},(N^{i,k_i})_{i\in\bbI,k_i\in\bbK^i}\big)$, 
{\ared with $(x,y)$ for the argument ``$x$''. }%
\end{example}

Denote by $\bbA$ a family of pairs of bounded measurable mappings $(a,b)$ on $\calx$ such that 
\beas
\sup_{(a,b)\in\bbA}\big(\big\||a|\big\|_\infty\vee\|b\|_\infty\big)\leq F
\eeas
for some {\ared positive} constant $F$. 
The true function $(a^*,b^*)$ correspond the true structure is assumed to satisfy $(a^*,b^*)\in\bbA$, 
as well as ${\mathfrak F}_T\subset\bbA$. 
We consider 
an estimator $(\wh{a}_T,\wh{b}_T)$ for $(a^*,b^*)$ from the data $(X_t,N_t)_{t\in[0,T]}$ for $T\in\bbT$ 
by optimizing $\Psi_T(a,b)$ {\ablue for a family ${\mathfrak F}_T$ of models} {\ared in $\bbA$}, e.g. deep learning {\ablue models}. 
The estimator $(\wh{a}_T,\wh{b}_T)$ takes the values in ${\mathfrak F}_T$. 

Let $(\ol{X},\ol{N})$ be an independent copy of $(X,N)$. 
The risk function {\ared (i.e., the expected prediction error)} when $(\wh{a}_T,\wh{b}_T)$ is used is 
\beas 
R_T
&=&
{\colorr E\bigg[T_1^{-1}\int_0^{T_1}\wh{U}_T(\ol{X}_t)dt\bigg]}
\eeas
for a fixed $T_1\in\bbT$, 
where 
\beas 
\wh{U}_T(x)&=&-\lambda^*(x)\cdot\big\{\wh{a}_T(x)-a^*(x)\big\}+\big\{\wh{b}_T(x)-b^*(x)\big\}. 
\eeas
We may choose $T_1={\tt h}$ due to the periodical stationarity. 
We also have the representation of $R_T$: 
\beas 
R_T
&=&
E\bigg[-T^{-1}\int_0^{T}\big\{{\ared\wh{a}_T}(\ol{X}_t)-a^*(\ol{X}_t)\big\}\cdot d\ol{N}_t+T^{-1}\int_0^{T}\big\{{\ared\wh{b_T}}(\ol{X}_t)-b^*(\ol{X}_t)\big\}dt\bigg]
\eeas
for $T\in\bbT$.

{\bred
The following compatibility condition is assumed: there exists a positive constant $\constzz\geq1$ such that 
{\bred 
\bea\label{0703240806}
\constzz^{-2}\big\{
\big|a(x)-a^*(x)\big|^2+\big|b(x)-b^*(x)\big|^2\big\}
&\leq&
-\lambda^*(x)\cdot\big\{a(x)-a^*(x)\big\}+\big\{b(x)-b^*(x)\big\}
\nn\\&\leq&
\constzz^2\big\{
\big|a(x)-a^*(x)\big|^2+\big|b(x)-b^*(x)\big|^2\big\}
\eea
for all $(a,b)\in{\bred\bbA}$, 
$T\in\bbT$, and all $x\in\calx$. %
Under (\ref{0703240806}), in particular, 
}
\bea\label{202412111052} &&
{\bred\constzz^{-1}}
\bigg({\tt h}^{-1}\int_0^{\tt h}E_{\ol{X}}\bigg[\big|a(\ol{X}_t)-a^*(\ol{X}_t)\big|^2
+\big|b(\ol{X}_t)-b^*(\ol{X}_t)\big|^2\bigg]dt\bigg)^{1/2}
\nn\\&\leq&
\bigg({\tt h}^{-1}\int_0^{\tt h}E_{\ol{X}}\bigg[-\lambda^*(\ol{X}_t)\cdot\big\{a(\ol{X}_t)-a^*(\ol{X}_t)\big\}
+\big\{b(\ol{X}_t)-b^*(\ol{X}_t)\big\}\bigg]dt\bigg)^{1/2}
\nn\\&\leq&
{\bred\constzz}
\bigg({\tt h}^{-1}\int_0^{\tt h}E_{\ol{X}}\bigg[\big|a(\ol{X}_t)-a^*(\ol{X}_t)\big|^2
+\big|b(\ol{X}_t)-b^*(\ol{X}_t)\big|^2\bigg]dt\bigg)^{1/2}
\eea
for all $(a,b)\in{\bred\bbA}$ 
{\ared and $T\in\bbT$}. 
{\ared Here $E_{\ol{X}}$ stands for the expectation with respect to ${\ol{X}}$.} %
Such a condition can be checked e.g. by the estimate: 
for any positive constants $x_0$ and $x_1$, there exist constants $c_0$ and $c_1$ such that 
\beas
c_0(x-1)^2\leq-\log x+x-1 \leq c_1(x-1)^2
\eeas
for all $x\in[x_0,x_1]$. 
}
In Example \ref{202412120701}, 
{\ablue when} the family of mappings $\lambda=(\lambda^i)_{i\in\bbI}$ associated with {\bred$(a.b)\in\bbA$} 
satisfies $0<\inf_{x\in\calx,T\in\bbT}\lambda^i(x)\leq\sup_{x\in\calx,T\in\bbT}\lambda^i(x)<\infty$, 
then the compatibility condition (\ref{0703240806}) holds true. 
The compatibility condition is a condition on the structure of {\bred$\bbA$}. 
{\bred The compatibility condition can be verified in a similar manner in Examples \ref{202504112106} and \ref{202504112107}. }%

Suppose that $\bbA$ admits a distance $\bbd$ such that 
\bea\label{202412251703}
\bbd\big((a',b'),(a,b)\big)
&\geq& 
2{\bred
\constzz}
{\ared\sfd_N}(1+\||\lambda^*|\|_\infty)(\|{\ared|}a'-a{\ared|}\|_\infty+\|b'-b\|_\infty)
\eea
{\ared for $(a,b), (a',b')\in\bbA$. }%
Then, in particular, 
\beas 
\bbd\big((a',b'),(a,b)\big)
&\geq& 
E\bigg[T_1^{-1}\int_0^{T_1}\big|a'(\ol{X}_t)-a(\ol{X}_t)\big||\lambda^*(\ol{X}_t)|dt+T_1^{-1}\int_0^{T_1}\big|b'(\ol{X}_t)-b(\ol{X}_t)\big|dt\bigg]. 
\eeas
Define $\Delta_T$ by 
\beas
\Delta_T 
&=& 
E\bigg[T^{-1}\Psi_T(\wh{a}_T,\wh{b}_T)-\inf_{(a,b)\in{\mathfrak F}_T}T^{-1}\Psi_T(a,b)\bigg]
\eeas

The $\alpha$-mixing coefficient for {\ared $\tth$-periodically stationary process} $X$ is given by 
{\ared 
\beas 
\alpha^X_\tth(k) &=& \sup_{j\in\bbZ_+}\sup_{A\in\calf_{[0,j\tth]}^X,B\in\calf_{[(j+k)\tth,\infty)}^X}\big|P[A\cap B]-P[A]P[B]\big|
\eeas
}%
for $k\in\bbZ_+$, 
where $\calf_I^X=\sigma[X_s;\>s\in I]$ for $I\subset\bbR_+$, i.e., 
the $\sigma$-field generated by $\{X_s;\>s\in I\}$. 
We assume that 
{\ared $\alpha^X_\tth(k)\leq {\ared\gamma}^{-1}e^{-{\ared\gamma}h}$ for all $k\in\bbZ_+$, }%
for some constant ${\ared\gamma}>0$. 
{\ared 
A usual $\alpha$-mixing coefficient for $X$ is}
\beas 
\alpha^X(h) &=& \sup_{t\in\bbR_+}\sup_{A\in\calf_{[0,t]}^X,B\in\calf_{[t+h,\infty)}^X}\big|P[A\cap B]-P[A]P[B]\big|
\eeas
{\ared 
If $\alpha^X(h)\leq{\ared\gamma'}^{-1}e^{-{\ared\gamma'}h}$ for all $h\in\bbR_+$, for some constant ${\ared\gamma'}>0$, then 
the $\alpha$-mixing coefficient $\alpha^X_\tth(k)$ geometrically decays. }

Let $\ttT=T/{\tt h}$ {\ared for $T\in\bbT$}. 
We give a rate of convergence of $R_T$. 
\begin{theorem}\label{0607170208}
{\bred Let $\xi$ be any positive number. Then }%
there exists a constant $\constq$ depending on ${\ared\gamma}$, {\bred $\tth$, $\big\||\lambda^*|\big\|_\infty$, $\sfd_N$, $\constzz$ and $\xi$}, such that 
\bea\label{202412220016}
R_T 
&\leq&
2\Delta_T+2\inf_{(a,b)\in{\mathfrak F}_T}{\tt h}^{-1}E\big[\Psi_{\tt h}(a,b)-\Psi_{\tt h}(a^*,b^*)\big]
+\constq(1+F^2)\bigg[{\bred\ttT}^{-1}(\log {\bred\ttT})^2\log\caln_T+\delta\bigg]
\nn\\&&
\eea
whenever 
{\bred $\ttT\geq2\vee\{\xi(\log\ttT)^2\log\caln_T\}$ and $\caln_T\geq2$. }%
Here {\ared$\caln_T=\caln_{T,\delta}$} is the covering number of ${\mathfrak F}_T$ 
by the $\delta$-balls with respect to the distance $\bbd$. 
\end{theorem}

We will prove Theorem \ref{0607170208} in Section \ref{0610080715}.

\section{Application to deep learning}\label{202412271344}
The inequality (\ref{202412220016}) can provide a rate of convergence of the risk 
if combined with an error bound of the approximation by the machine ${\mathfrak F}_T$ and 
an estimate of its covering number $\caln_T$. 
Schmidt-Hieber \cite{schmidt2020nonparametric} considered a deep neural network 
with ReLU activation function and presented a covering number when the network is fitted under a sparse condition. 

The shifted ReLU activation function $\sigma_v:\bbR^\sfd\to\bbR^\sfd$ is defined as 
\beas 
\sigma_v(x) &=& \big((x_1-v_1)_+,...,(x_\sfd-v_\sfd)_+\big)^\star
\eeas
where for $x=(x_1,...,x_\sfd)^\star\in\bbR^\sfd$ and ${\bf v}=(v_1,...,v_\sfd)^\star$,  $x_+\max\{u,0\}$ for $u\in\bbR$. 
For weight matrices $W_i\in\bbR^{\sfp_{i+1}}\otimes\bbR^{\sfp_i}$ ($i=0,1,...,L$) and 
shift vectors ${\bf v}_i\in\bbR^{\sfp_i}$ ($i=1,...,L$), 
the mapping $f\big(\cdot;(W_L,...,W_1,W_0),({\bf v}_L,...,{\bf v}_1)\big):\bbR^{\sfp_0}\to\bbR^{\sfp_{L+1}}$ is defined as 
\bea\label{202412251346}
f\big(x;(W_L,...,W_1,W_0),({\bf v}_L,...,{\bf v}_1)\big)
&=&
W_L\sigma_{{\bf v}_L}W_{L-1}\sigma_{{\bf v}_{L-1}}\cdots W_1\sigma_{{\bf v}_1}W_0x
\quad(x\in\bbR^{\sfp_0}). 
\eea
The dimension $\sfp_0=\sfd_X$ of the input process $X$,  and $\sfp_L=1$ in the applications to point processes in this article
The set of functions $f\big(x;(W_L,...,W_1,W_0),({\bf v}_L,...,{\bf v}_1)\big) $ taking the form (\ref{202412251346}) is denoted by $\cald$, 
{\ablue and it is called a deep neural network or deep learning}. 
{\ablue Some restrictions are posed on $\cald$, 
depending on the situation one is working. }%
Schmidt-Hieber \cite{schmidt2020nonparametric} uses the class 
\beas 
{\mathfrak F}_T
&=& 
\bigg\{f\in\cald\text{ of the form }(\ref{202412251346});\>
\max_{\ell=0,...,L, j=1,...,L}\big(\|W_\ell\|_\infty\vee\|{\bf v}_j\|_\infty\big)\leq1,
\nn\\&&
\sum_{\ell=0}^L\|W_\ell\|_0+\sum_{j=1}^L\|{\bf v}_j\|_0\leq s,\>\|f\|_\infty\leq F
\bigg\}
\eeas
for a given positive constant $F$, where the parameters $L$ and $\sfp_i$ ($i=1,...,L$) 
determining the size of the learning machine, 
as well as the sparsity index $s$, depend on $T$. 
The $0$-norm $\|\cdot\|_0$ denotes the number of non-zero entries of the object. 

The function $g$ that generates the data is assumed to be expressed as a composite of functions of 
H\"older classes, in Schmidt-Hieber \cite{schmidt2020nonparametric}. 
Therein prepared is a ball of $\beta$-H\"older functions with radius $K$ denoted by 
\beas 
C^\beta(D,K)
&=&
\bigg\{{\tt g}:D\to\bbR;\>
\sum_{\alpha:|\alpha|<\beta}\|\partial^\alpha{\tt g}\|_\infty+
\sum_{\alpha:|\alpha|=\lfloor\beta\rfloor}
\sup_{x,y\in D,x\not=y}\frac{|\partial^\alpha{\tt g}(x)-\partial^\alpha{\tt g}(y)|}{\|x-y\|_\infty^{\beta-\lfloor\beta\rfloor}}\leq K\bigg\}
\eeas
for 
a domain $D$ in a Euclidean space and a number $K$, where
$\lfloor\beta\rfloor$ denotes the largest integer strictly smaller than $\beta$. 
A family $\calg=\calg(q,{\bf d},{\bf t},{\bm \beta},K)$ of the possible data generating mechanisms is a collection of functions $g={\tt g}_q\circ\cdots\circ{\tt g}_0$ 
such that 
${\tt g}_i=({\tt g}_{ij})_{j}:[a_i,b_i]^{d_i}\to[a_{i+1},b_{i+1}]^{d_{i+1}}$, where each component ${\tt g}_{ij}$ is a function of 
some of the arguments in $\bbR^{d_i}$ and satisfies 
${\tt g}_{ij}\in C^{\beta_i}([a_i,b_i]^{t_i},K)$, 
given vectors ${\bf d}=(d_0,...,d_{q+1})$, ${\bf t}=(t_0,...,t_q)$ and ${\bm\beta}=(\beta_0,...,\beta_q)$. 

In order to obtain a good bound for the risk $R_T$, a set of conditions is required to make ${\mathfrak F}_T$ 
sufficiently rich and not too large in the same time. 
Naturally, such conditions involve the smoothness of the target function $g$. 
As Schmidt-Hieber \cite{schmidt2020nonparametric}, we impose the following conditions:
\bea\label{202412251723}&&
F\geq\max\{K,1\}, \quad
\sum_{i=0}^q\log_2 4(t_i\vee\beta_i)\log_2T\leq L\simleq T\phi_T, \quad
\nn\\&&
T\phi_T\simleq\min\{\sfp_1,...,\sfp_L\}, \quad
s\asymp T\phi_T\log T. 
\eea

The effective smoothness {\tred index} is defined as 
$\beta^*_i=\beta_i\prod_{j=i+1}^q(\beta_j\wedge1)$, and the key convergence rate as 
\beas 
\phi_T &=& \max_{i=0,...,q}T^{-\frac{2\beta_i^*}{2\beta_i^*+t_i}}. 
\eeas

As Inequality (26) of Schmidt-Hieber \cite{schmidt2020nonparametric}, we obtain 
\bea\label{202412251655}
\inf_{(a,b)\in{\mathfrak F}_T}{\tt h}^{-1}E\big[\Psi_{\tt h}(a,b)-\Psi_{\tt h}(a^*,b^*)\big]
&\simleq&
\phi_T
\eea
by the compatibility (\ref{202412111052}). 
On the other hand, Lemma 5 of Schmidt-Hieber \cite{schmidt2020nonparametric} 
gives an estimate of the covering number as 
\bea\label{202412251659}
\log\caln_T
&\leq&
(s+1)\log\bigg[2\delta^{-1}(L+1)\prod_{\ell=0}^{L+1}(\sfp_\ell+1)\bigg]. 
\eea
The above covering number is based on the uniform norm but we can now take a metric $\bbd$ of (\ref{202412251703}) 
compatible with the sup-norm.

Following  Schmidt-Hieber \cite{schmidt2020nonparametric}, we obtain the following estimate of the risk in 
the prediction with ${\mathfrak F}_T$ specified above, if combined with 
the properties (\ref{202412251655})-(\ref{202412251659})  
($\sfp_\ell$ are bounded by $s$). 
\begin{theorem}\label{202412271233}
{\bred Let $\xi>0$. }%
If $\Delta_T\leq C_0\phi_TL(\log T)^4$ $(T\geq T_0)$ for some positive constants $C_0$ and $T_0>1$, then there exists a constant $C$ such that 
\bea\label{202412271236}
R_T &\leq& C\phi_T L(\log T)^4
\eea
for $T\geq T_0$ 
{\bred whenever $T\geq \xi(\log T)^2(s+1)\log\big[2\delta^{-1}(L+1)\prod_{\ell=0}^{L+1}(\sfp_\ell+1)\big]$.}
\end{theorem}
%

\begin{remark}\rm
\bd
\im[(i)] 
The error bound (\ref{202412271236}) has the factor $(\log T)^4$ instead of $(\log T)^2$ in Schmidt-Hieber\cite{schmidt2020nonparametric}. 
This factor comes from the large deviation estimate for a functional in the mixing condition. 
{\bred 
The error bound (\ref{202412271236}) becomes $C\phi_T(\log T)^5$ when $L\asymp\log T$. 
}
\im[(ii)] The error bound (\ref{202412271236}) is minimax-optimal up to the logarithmic factor in that 
our model includes the case where 
the covariate process $X$ is periodically independent. 
Schmidt-Hieber\cite{schmidt2020nonparametric} showed the optimality for an independent input process in the context of nonparametric {\tred regression} setting and the risk function is the same under the compatibility condition. 
\im[(iii)] 
Suzuki and Nitanda \cite{suzuki2021deep} propose the use of an anisotropic Besov space to represent the target function. 
This approach can be taken to obtain a better error bound also for the point process models. 
\ed
\end{remark}

\section{The marked ratio model}\label{202412271346}

\subsection{A simulation study}
We propose in this section a simulation study in the case of marked intensities. Recall that in this case we consider the marked intensity processes
\begin{equation}
	\lambda^{i,k_i}(X_t,Y_t) = \lambda_0(t) \lambda^i(X_t)p_i^{k_i}(Y_t), 
	\quad 
	i\in\mathbb I, k_i\in\mathbb K_i,
\label{eq:simulation_model_eq1}
\end{equation}
where we assume that $\sum_{k_i\in\mathbb K_i} p_i^{k_i} = 1$ for all $i\in\mathbb I$, 
{\ablue and $\lambda_0$ is an unobserved baseline intensity. 
We refer the reader to Muni Toke and Yoshida \cite{muni2022marked} for more details on the model. }

\subsubsection{Description of the numerical example}
\label{subsec:SimulationStudy-Model}

In our numerical example we consider a 4-dimensional process with 2-dimensional marks, i.e. we set $\sfd_N=4$, $\mathbb I=\{0,1,2,3\}$ and for all $i\in\mathbb I$,  $\mathbb K_i=\{0,1\}$. Covariates process $X$ is $\sfd_X$-dimensional with $\sfd_X=2$, covariates process $Y$ is $\sfd_Y$-dimensional with $\sfd_Y=1$, and all three coordinate covariates are independent Ornstein-Uhlenbeck (OU) processes. More precisely, we consider $(B^{X^0}, B^{X^1}, B^{Y})$ a 3-dimensional Brownian motion in our probability space and set:
\begin{equation}\begin{cases}
	dX^j_t = \theta_{X^j} (\bar x_{X^j}-X^j_t)\,dt + \sigma_{X^j}\,dB^{X^j}_t, & j=0,1,
	\\
	dY_t = \theta_{Y} (\bar x_{Y}-Y_t)\,dt + \sigma_{Y}\,dB^{Y}_t.	
	&
	\end{cases}
\end{equation}
Values of the OU parameters $((\theta_{X^j},\bar x_{X^j},\sigma_{X^j})_{j=0,1}, \theta_{Y}, \bar x_{Y},\sigma_{Y})$ are given in Table \ref{table:OUparameters}.
\begin{table}
\begin{center}
\begin{tabular}{|c|ccc|}
	\hline
	& $\theta_{\cdot}$ & $\bar x_{\cdot}$ & $\sigma_{\cdot}$
	\\ \hline
	$X^0$ & 0.1 & 0.0 & 0.1
	\\
	$X^1$ & 0.2 & 0.0 & 0.2
	\\
	$Y$ &  0.1 & 0.0 & 0.1
	\\ \hline
\end{tabular}
\end{center}
\caption{Numerical values for the OU covariate processes.}
\label{table:OUparameters}
\end{table}
We keep the number of covariates reasonably low in this numerical example so that our fitting results can still be represented graphically in a manageable way.

{\bred The base line intensity is set to $\lambda_0 (t) = 1 + \cos ( 2 \pi t )$. }%
Non-marked intensities $\lambda^i$ are defined for $x=(x_0,x_1)$ as:
\begin{equation}\begin{cases}
	\lambda^0(x) = 2 + \tanh(x_0)\exp(-x_1^2),
	\\
	\lambda^1(x) = 2 + \cos(\pi x_0) \tanh(x_1),
	\\
	\lambda^2(x) = 2 + \sin(2\pi x_0) e^{x_1}(1+e^{x_1})^{-1},
	\\
	\lambda^3(x) = 3 - \exp(-x_0^2),
\end{cases}
\end{equation}
and mark probabilities $p_i^{k_i}$ are written:
\begin{equation}\begin{cases}
	p_0^0(y) = 0.25, & p_0^1 = 1-p_0^0,
	\\
	p_1^0(y) = 0.05 + 0.9 |\cos(\pi y)|, & p_1^1 = 1-p_1^0,
	\\
	p_2^0(y) = e^{y}(1+e^{y})^{-1}, & p_2^1 = 1-p_2^0,
	\\
	p_3^0(y) = 0.6 \exp(-y^2), & p_3^1 = 1-p_3^0.
\end{cases}
\label{eq:simulation_model_eq4}
\end{equation}
Numerical simulations of the point processes $(N^{i,k_i})_{i\in\mathbb I, k_i\in\mathbb K_i}$ are carried via a thinning algorithm.

\subsubsection{One-step ratio estimation}
\label{subsec:oneStepEstimation}

We define a first estimation method in the spirit of \cite{muni2019analyzing}. We start by considering the $8$-dimensional point process $(N^{i,k_i})_{i\in\mathbb I, k_i\in\mathbb K_i}$. We define the ratio functions
\begin{equation}
	r_1^{i,k_i}(x,y) = \frac{\lambda^{i,k_i}(x,y)}{\sum_{j\in\mathbb I, k_j\in\mathbb K_j} \lambda^{j,k_j}(x,y)} 
	\quad\text{ and }\quad
	\tilde r_1^{i,k_i}(x,y) = \frac{r_1^{i,k_i}(x,y)}{r_1^{0,0}(x,y)}.
\end{equation}
Obviously, $\sum_{i\in\mathbb I, k_i\in\mathbb K_i}r_1^{i,k_i}=1$, $\tilde r_1^{0,0}=1$ and $\sum_{i\in\mathbb I, k_i\in\mathbb K_i} \tilde r_1^{i,k_i}=\frac{1}{r_1^{0,0}}$.
In this first estimation method, we set
\begin{equation}
	l_1^{i,k_i}(x,y) = \log \tilde r_1^{i,k_i}(x,y)
\end{equation}
and these functions are estimated for $(i,k_i)\neq(0,0)$ with a neural network.

We define a standard dense feed-forward neural network with a $(d_X,n_1^N)$-shaped input layer for the covariates $X$, $n^L_1$ inner layers with $n_1^N$ neurons per layer and a final $(n^N_1,7)$-shaped output layer to output the estimated quantities $\hat l_1^{i,k_i}(x,y)$, $(i,k_i)\neq(0,0)$, that approximate the $l_1^{i,k_i}(x,y)$. All layers except the last hidden one and the output one use a LeakyReLu activation function.
In the general terminology of the previous sections, the contrast function $\Psi_T(a_1,b_1)$ is in this case defined with $b_1(x,y)=0$ and
\begin{align}
	a_1^{i,k_i}(x,y) 
	= \log r_1^{i,k_i}(x,y)
	& = \log \frac{\tilde r_1^{i,k_i}(x,y)}{\sum_{j\in\mathbb I, k_j\in\mathbb K_j} \tilde r_1^{j,k_j}(x,y)}.
\end{align}
The loss function $\mathcal L_1$ of the neural network computed on a sample $\mathcal S_{1,T}=\{(X_t,Y_t,(N^{i,k_i}_t)_{i,k_i})\}_{t\in[0,T]}$ is thus
\begin{equation}
	\mathcal L_1(\mathcal S_{1,T})
	=
	- \int_0^T
	\sum_{
		i\in\mathbb I, k_i\in\mathbb K_i
	}
	\log \frac{\exp(l_1^{i,k_i}(X_t,Y_t))}{\sum_{j\in\mathbb I, k_j\in\mathbb K_j} \exp(l_1^{j,k_j}(X_t,Y_t))} dN^{i,k_i}_t.
\end{equation}
Recall that $l_1^{0,0}=0$ is not learned. Index $1$ in the notation of this section is used to indicate that this is our first estimation method (one-step ratio estimation).

\subsubsection{Two-step ratio estimation}
\label{subsec:twoStepEstimation}

We now define a second estimation method in the spirit of \cite{muni2022marked}. In a first step we use a ratio model on the non-marked intensities $\lambda^i(X_t)$ and then in a second step we use $\#\mathbb I=4$ other ratio estimations on the mark probabilities $p_i^{k_i}$, one for each $i\in\mathbb I$. 

The notation for the first step of this second estimation is 
\begin{equation}
	r_2^{i}(x) = \frac{\lambda^{i}(x)}{\sum_{j\in\mathbb I} \lambda^{j}(x)},
	\quad
	\tilde r_2^{i}(x) = \frac{r_2^{i}(x)}{r_2^{0}(x)}
	\quad\text{ and }\quad
	l_2^{i}(x) = \log \tilde r_2^{i}(x)
\end{equation}
and these functions are estimated for $i\neq0$ with a neural network. In order to compute the estimators $\hat l_2^i(x)$ we use the general architecture previously defined in the first estimation method, but now with parameters $n^L_{2,1}, n^N_{2,1}$ and a 3-dimensional output ($i\in\mathbb I\setminus\{0\}$).
The loss function $\mathcal L_{2,1}$ of the neural network computed on a sample $\mathcal S_{2,1,T}=\{(X_t,(N^{i}_t)_{i})\}_{t\in[0,T]}$ is thus
\begin{equation}
	\mathcal L_{2,1}(\mathcal S_{2,1,T})
	=
	- \int_0^T
	\sum_{i\in\mathbb I}
	\log \frac{\exp(l_2^{i}(X_t))}{\sum_{j\in\mathbb I} \exp(l_2^{j}(X_t))} dN^{i}_t.
\end{equation}

The notation for the $i$-th ratio model of the second step of the second estimation method is then
\begin{equation}
	r_2^{i,k_i}(y) = \frac{p_i^{k_i}(y)}{\sum_{k\in\mathbb K_i} p_i^{k}(y)} = p_i^{k_i}(y),
	\quad
	\tilde r_2^{i,k_i}(y) = \frac{r_2^{i,k_i}(y)}{r_2^{i,0}(y)}
	= \frac{p_i^{k_i}(y)}{p_i^{0}(y)}
	\quad\text{ and }\quad
	l_2^{i,k_i}(y) = \log \tilde r_2^{i,k_i}(y)
\end{equation}
and these functions are estimated for $k_i\neq0$ with a neural network. Again, we use in order to compute the estimators $\hat l_2^{i,k_i}(y)$ the same general architecture for the neural networks, now with a $d_Y$-dimensional input, parameters $n^L_{2,2}, n^N_{2,2}$ and a 1-dimensional output ($k_i\in\mathbb K_i\setminus\{0\}$).
The loss function $\mathcal L^i_{2,2}$ of the neural network computed on a sample $\mathcal S^i_{2,2,T}=\{(Y_t,(N^{i,k_i}_t)_{k_i})\}_{t\in[0,T]}$ is in this case
\begin{equation}
	\mathcal L^i_{2,2}(\mathcal S^i_{2,2,T})
	=
	- \int_0^T
	\sum_{k_i\in\mathbb K}
	\log \frac{\exp(l_2^{i,k_i}(Y_t))}{\sum_{k\in\mathbb K} \exp(l_2^{i,k}(Y_t))} dN^{i,k_i}_t.
\end{equation}

\subsubsection{Fitting results}
\label{subsec:SimulationStudy-FittingResults}
As a first illustration, we simulate the model \eqref{eq:simulation_model_eq1}-\eqref{eq:simulation_model_eq4} for an horizon $T=128,000$ (note that given the above definitions, a sample has roughly $2T$ points in each of the 4 dimensions of the process in this model). We then fit the model with our two estimation methods and the parameters $n^L_1=n^L_{2,1}=n^L_{2,2}=8$ and $n^N_1=n^N_{2,1}=n^LN_{2,2}=64$. 

\begin{figure}
	\includegraphics[width=\textwidth]{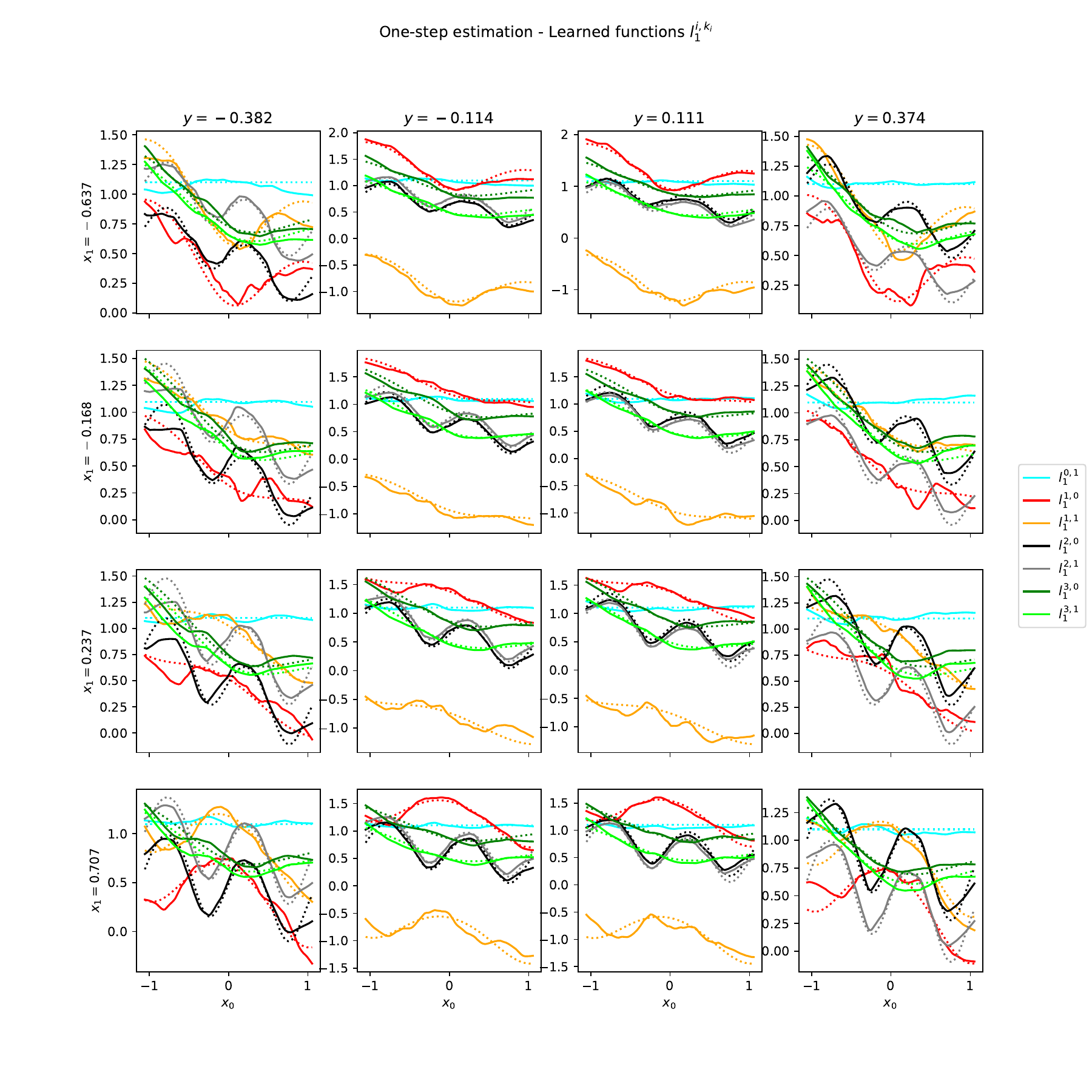}
	\caption{Simulation study --- Estimated functions $\hat l_1^{i,k_i}(x,y)$ by the one-step estimation method. True functions are plotted as dotted lines of the color of the corresponding estimated function.}
	\label{fig:OneStepEstimation-LearnedFunctions-128000}
\end{figure}

\begin{figure}
	\includegraphics[width=\textwidth]{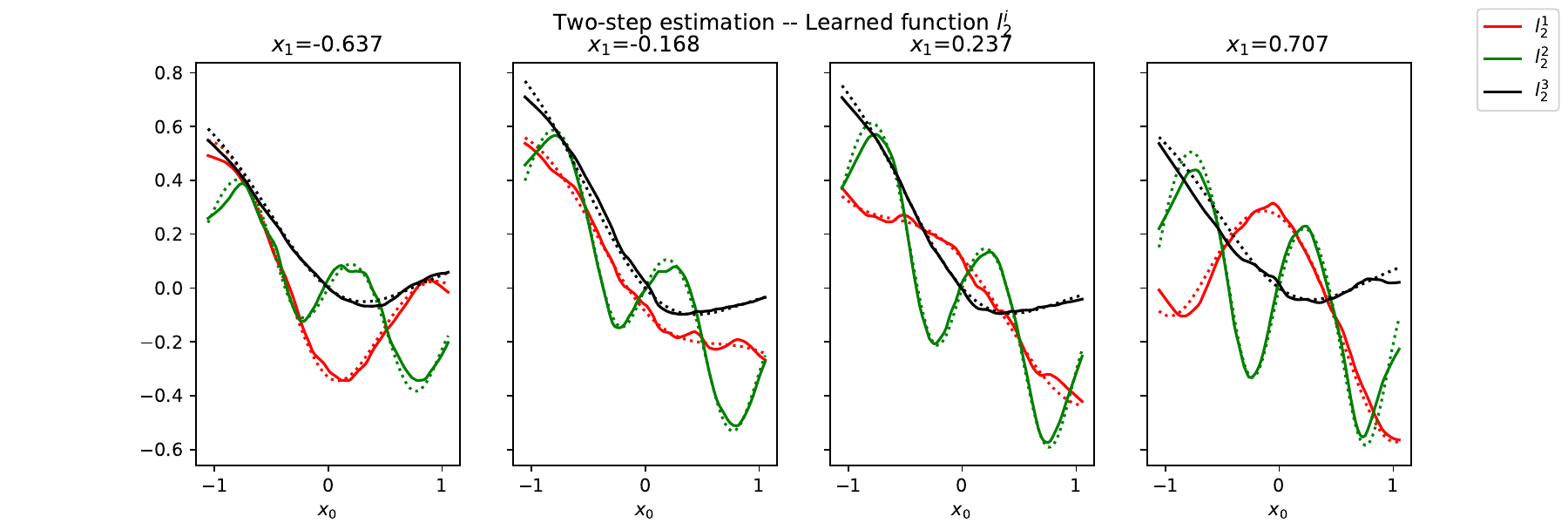}
	\caption{Simulation study --- Estimated functions $\hat l_2^{i}(x)$ by the two-step estimation method. True functions are plotted as dotted lines of the color of the corresponding estimated function.}
	\label{fig:TwoStepEstimation-LearnedFunctions-128000}
\end{figure}

\begin{figure}
	\includegraphics[width=\textwidth]{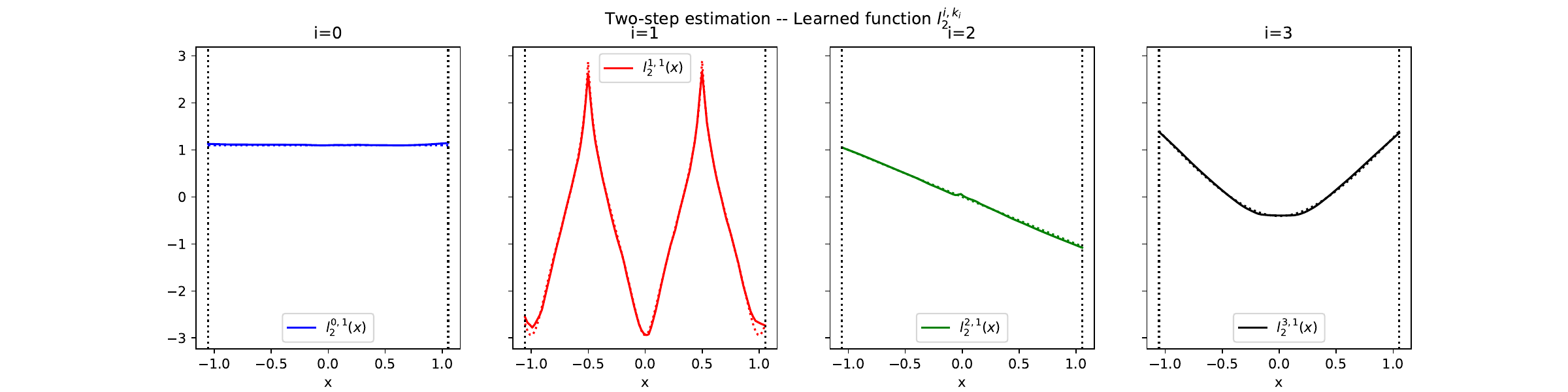}
	\caption{Simulation study --- Estimated functions $\hat l_2^{i,k_i}(y)$ by the two-step estimation method. True functions are plotted as dotted lines of the color of the corresponding estimated function.}
	\label{fig:TwoStepEstimation-LearnedFunctionsProba-128000}
\end{figure}

Figure \ref{fig:OneStepEstimation-LearnedFunctions-128000} plots the true functions $l_1^{i,k_i}(x,y)$ and the estimated functions $\hat l_1^{i,k_i}(x,y)$ by the one-step estimation method. In order to better visualize the results, we provide plots of the 7 functions $x_0\mapsto \hat l_1^{i,k_i}(x_0, \hat q_{X^1}(\alpha), \hat q_{Y}(\beta))$ where $\hat q_{X^1}(\alpha)$ is the $\alpha$-quantile of the empirical distribution of $X^1$, $\hat q_{Y}(\beta))$ is the $\beta$-quantile of the empirical distribution of $Y$, and $\alpha,\beta\in [0.2, 0.4, 0.6, 0.8]$, hence the $4\times 4$ matrix of plots. 
Figure \ref{fig:TwoStepEstimation-LearnedFunctions-128000} plots the true functions $l_2^{i}(x)$ and the estimated functions $\hat l_2^{i}(x)$ (Again, we plot $x_0\mapsto \hat l_2^{i}(x_0, \hat q_{X^1}(\alpha))$ for $\alpha\in [0.2, 0.4, 0.6, 0.8]$). Figure \ref{fig:TwoStepEstimation-LearnedFunctionsProba-128000} plots the true functions $l_2^{i,k_i}(y)$ and the estimated functions $\hat l_2^{i,k_i}(y)$. 
All {\ared these} graphs illustrate the ability of the estimation methods to retrieve the various shapes of ratio functions defined by the model.

Now, in order to compare the estimation methods, we illustrate the results in terms of probabilities. In the model \eqref{eq:simulation_model_eq1}-\eqref{eq:simulation_model_eq4}, the probability that an observed event in state $(x,y)$ is of type $i$ and mark $k_i$ is
\begin{equation}
	p^{i,k_i}(x,y)=\frac{\lambda^{i,k_i}(x,y)}{\sum_{j\in\mathbb I, k_j\in\mathbb K_j}\lambda^{j,k_j}(x,y)}.
\label{eq:piki_joint}
\end{equation}
Note that $p^{i,k_i}$ and $p_{i}^{k_i}$ are not the same. $p^{i,k_i}$ defined at Equation \eqref{eq:piki_joint} {\ablue is the joint probability of the type $i$ and the mark $k_i$}, while $p_{i}^{k_i}$ defined at Equation \eqref{eq:simulation_model_eq1} {\ablue is the conditional probability of the mark $k_i$ given the type $i$}.

\begin{figure}
	\includegraphics[width=\textwidth]{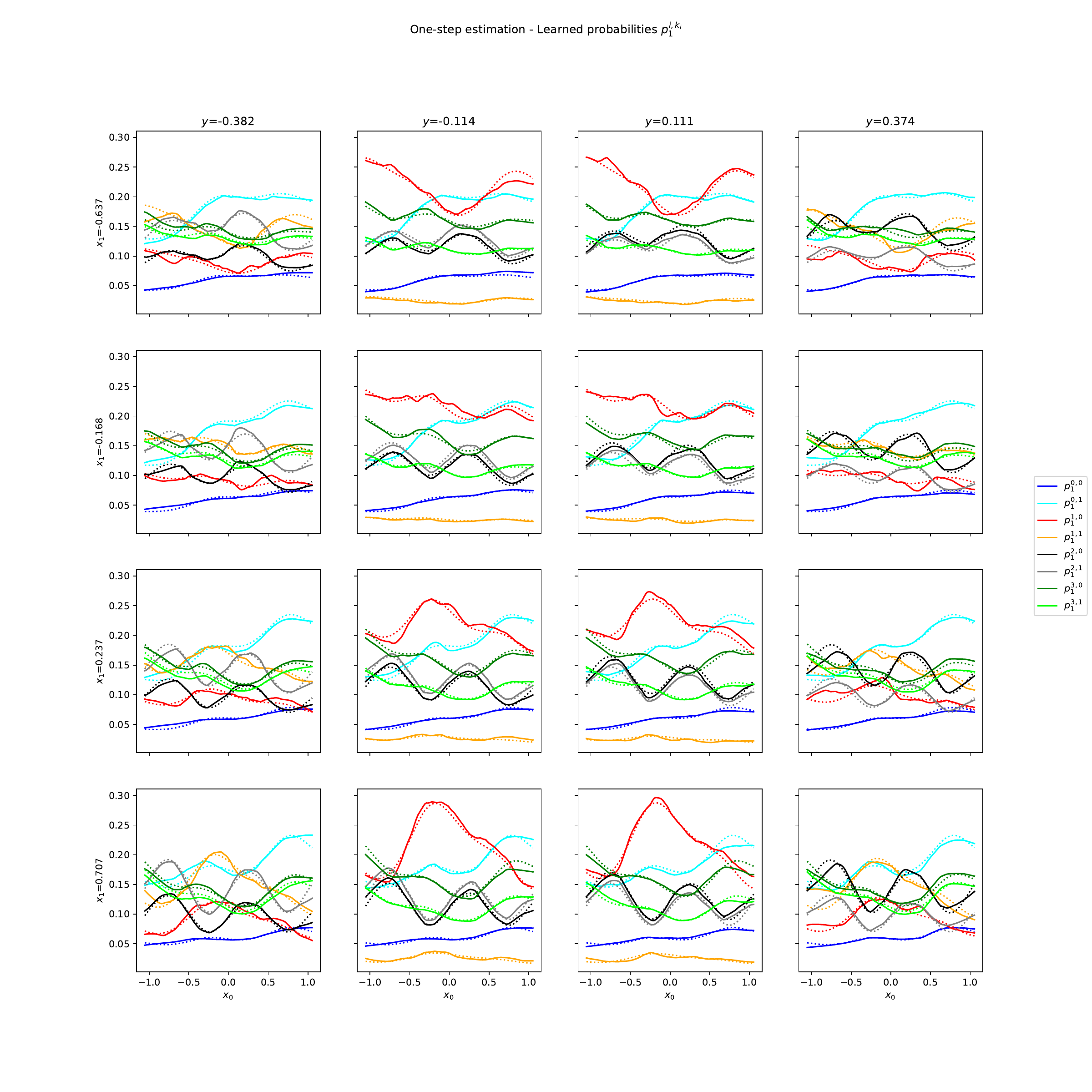}
	\caption{Simulation study --- Estimated probabilities $\hat p_1^{i,k_i}(x,y)$ by the one-step estimation method. True functions are plotted as dotted lines of the color of the corresponding estimated function.}
	\label{fig:OneStepEstimation-LearnedProbabilities-128000}
\end{figure}

\begin{figure}
	\includegraphics[width=\textwidth]{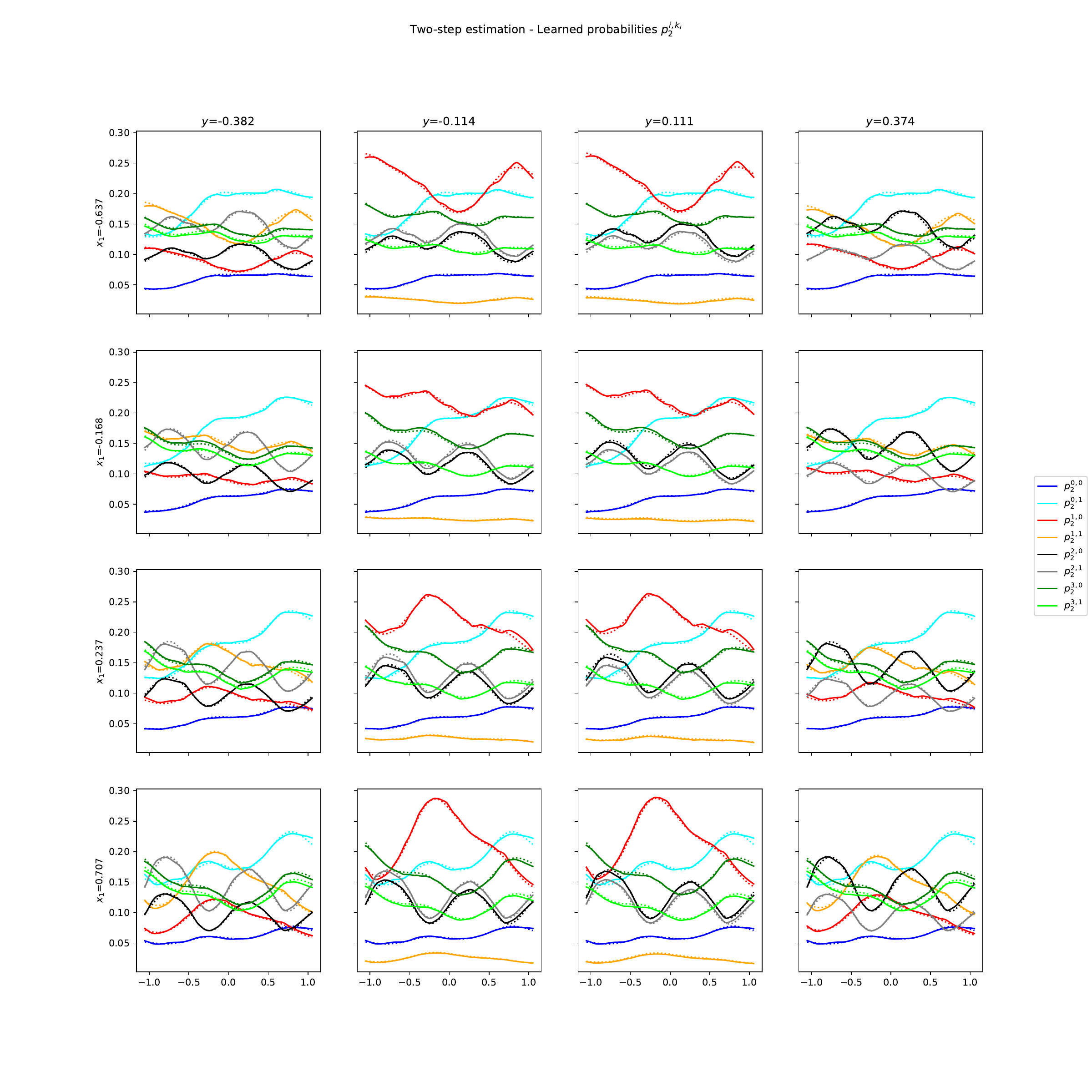}
	\caption{Simulation study --- Estimated probabilities $\hat p_2^{i,k_i}(x,y)$ by the two-step estimation method. True functions are plotted as dotted lines of the color of the corresponding estimated function.}
	\label{fig:TwoStepEstimation-LearnedProbabilities-128000}
\end{figure}

Probabilities $p^{i,k_i}(x,y)$ are straightforwardly estimated by the one-step estimation method with
\begin{equation}
	\hat p_1^{i,k_i}(x,y)=\frac{\exp(\hat l_1^{i,k_i}(x,y))}{\sum_{j\in\mathbb I, k_j\in\mathbb K_j}\exp(\hat l_1^{j,k_j}(x,y))},
\end{equation}
and by the two-step estimation method with
\begin{equation}
	\hat p_2^{i,k_i}(x,y)=
	\frac{\exp(\hat l_2^{i}(x))}{\sum_{j\in\mathbb I}\exp(\hat l_2^{j}(x))}	
	\frac{\exp(\hat l_2^{i,k_i}(y))}{\sum_{k\in\mathbb K_i}\exp(\hat l_2^{i,k}(y))}.
\end{equation}
Figure \ref{fig:OneStepEstimation-LearnedProbabilities-128000} plots the true functions $p^{i,k_i}(x,y)$ and the estimated functions $\hat p_1^{i,k_i}(x,y)$ and Figure \ref{fig:TwoStepEstimation-LearnedProbabilities-128000} plots the true functions $p^{i,k_i}(x,y)$ and the estimated probabilities $\hat p_2^{i,k_i}(x,y)$. We use the  $4\times 4$-matrix representation defined above for Figure \ref{fig:OneStepEstimation-LearnedFunctions-128000}.
Both methods provide visually high-quality fits for the event probabilities of the model. However, a careful examination of the plots indicates that the two-step estimation method estimates provide a better fit {\ared to} the true probabilities. We formalize this observation in the following section.

\subsubsection{Comparison of the methods and convergence results}
\label{subsec:SimulationStudy-Convergence}

Recall that using the general terminology of the Section \ref{202412271342} the risk function of our models is (since $b\equiv 0$ in the ratio estimations)
\begin{equation}
	R_T = E\left[ - \frac{1}{T} \int_0^ T
		\left\{ \hat a(\ol X_t,\ol Y_t)
			- a^*(\ol X_t,\ol Y_t)
		\right\} \cdot d\ol N_t,
	\right]
\end{equation}
where $(\ol X, \ol Y, \ol N)$ are independent copies of $(X,Y,N)$. We can thus simulate a new sample of length $T$ of our model and compute empirical versions $\mathcal R_T$ of the risk functions.

In the one-step estimation, we obtain on the sample $\ol{\mathcal S}_{1,T}=\{(\ol X_t,\ol Y_t,(\ol N^{i,k_i}_t)_{i,k_i})\}_{t\in[0,T]}$
\begin{align}
	\mathcal R_{1,T}(\ol{\mathcal S}_{1,T}) & = 
	- \frac{1}{T} \int_0^T
	\left\{ \log \hat r_1(\ol X_t,\ol Y_t)
	- \log r_1(\ol X_t,\ol Y_t) \right\} d\ol N_t
	\nonumber
	\\ & = - \frac{1}{T} \int_0^T
		\sum_{
		i\in\mathbb I, k_i\in\mathbb K_i
	}
	\log \frac{\exp(\hat l_1^{i,k_i}(\ol X_t,\ol Y_t))}{\sum_{j\in\mathbb I, k_j\in\mathbb K_j} \exp(\hat l_1^{j,k_j}(\ol X_t,\ol Y_t))} 
	\frac{\sum_{j\in\mathbb I, k_j\in\mathbb K_j} \exp(l_1^{j,k_j}(\ol X_t,\ol Y_t))}{\exp(l_1^{i,k_i}(\ol X_t,\ol Y_t))}
	d\ol N^{i,k_i}_t.
\end{align}

The empirical risk functions $\mathcal R_{2,1,T}(\ol{\mathcal S}_{2,1,T})$ and $\mathcal R^i_{2,2,T}(\ol{\mathcal S}^i_{2,2,T})$, $i\in\mathbb I$ are analogously defined with appropriate subsamples $\ol{\mathcal S}^i_{2,1,T}$ and $\ol{\mathcal S}^i_{2,2,T}$.

Moreover, to provide a complementary view, we define a standard uniform mean square error $\epsilon_{L^2,m}$ of the estimation method $m=1,2$ ($m=1$ for the one-step ratio method, $m=2$ for the two-step  ratio method). For each covariate $Z\in\{X^0, X^1,Y\}$, we compute the $1\%$ and $99\%$ empirical quantiles $q^Z_{0.01}$ and $q^Z_{0.99}$ on the (full) data and define a 1-dimensional regular grid of size $G+1$:
\begin{equation}
	\mathcal G^Z = \left\{
		q^Z_{0.01} + g \frac{q^Z_{0.99}-q^Z_{0.01}}{M} : g=0,\ldots,G
	\right\}
\end{equation}
The uniform $L^2$-type error is straightforwardly defined on the 3-dimensional grid $\mathcal G = \mathcal G^{X_0}\times \mathcal G^{X_1}\times \mathcal G^{Y}$ as
\begin{equation}
	\epsilon_{L^2,m} = \sum_{i\in\mathbb I} 
	\sum_{k_i\in\mathbb K_i}
	\sqrt{ \frac{1}{\#\mathcal G}
		\sum_{(x,y)\in\mathcal G}
		\left( \hat p_m^{i,k_i}(x,y) - p^{i,k_i}(x,y) \right )^2
	}
\end{equation}
This error is uniform in the sense that it does not take into account the distribution of the covariates. Similarly, a $L^\infty$-type error on the regular grid $\mathcal G$ is defined as
\begin{equation}
	\epsilon_{L^\infty,m} = \max_{i\in\mathbb I} 
	\max_{k_i\in\mathbb K_i}\max_{(x,y)\in\mathcal G} 
	\left| \hat p_m^{i,k_i}(x,y) - p^{i,k_i}(x,y) \right|.
\end{equation}

\begin{figure}
	\includegraphics[width=\textwidth]{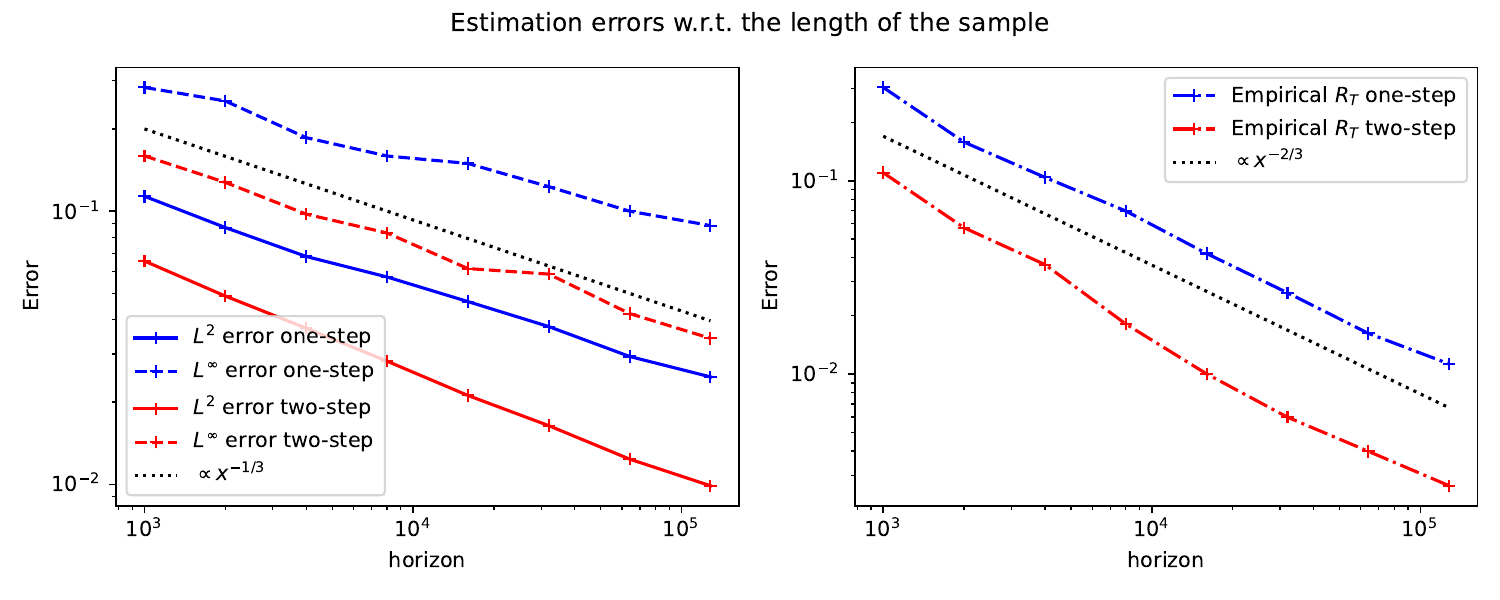}
	\caption{Simulation study --- $L^2$-errors (full lines, left panel), $L^\infty$-errors (dashed lines, left panel) and empirical risk function $\mathcal R_T$ (dash-dotted lines, right panel) as function of the horizon of the simulation for the one-step (blue) and the two-step (red) estimation methods. Dotted black lines with slopes $-1/3$ and $-2/3$ are plotted for visual guidance.}
	\label{fig:Horizon-Analysis}
\end{figure}

Figure \ref{fig:Horizon-Analysis} plots these three measures of estimation error as function of the simulation horizon in the case of the one-step and the two-step estimation methods. For each horizon $T$, we simulate 20 samples and run both estimation methods on each sample. We then compute the mean $L^2$-errors, mean $L^\infty$-errors and mean empirical risk function $\mathcal R_T$ across the 20 estimations.
Both methods exhibit close order of convergence with respect to the length of the sample, which is close to $-1/3$ for $L^2$ and $L^\infty$ errors and $-2/3$ in the case of the risk function $\mathcal R_T$. The superiority of the two-step estimation method, which takes into account the multiplicative structure of the model, is clear.

\subsubsection{Robustness with respect to shapes of the neural networks}
\label{subsec:SimulationStudy-Robustness}

Results of the previous sections have been obtained with the parameters $n^L_1=n^L_{2,1}=n^L_{2,2}=8$ and $n^N_1=n^N_{2,1}=n^N_{2,2}=64$.
We now run some tests to illustrate the robustness of the estimation with respect to the architecture of the neural networks used. For number of inner layers $n^L \in\{1,2,4,6,8,10, 12, 16, 20\}$ and each number of neurons per layer $n^N\in\{4, 8, 16, 32, 64, 128, 256, 512\}$, we run 20 simulations with horizon $T=32,000$ and estimations using both estimation methods. In the two-step estimations methods, all networks use the same parameters, i.e. $n^L_{2,1}=n^L_{2,2}=n^L$ and $n^N_1=n^N_{2,1}=n^N$.
Figures \ref{fig:Tuning-L2}, \ref{fig:Tuning-Linf} and \ref{fig:Tuning-RT} in Appendix provide the robustness results with respect to the shapes of the neural networks as heatmaps of the three error measures defined in Section \ref{subsec:SimulationStudy-Convergence}. It appears clearly that the estimation methods are quite robust with respect to the shapes of the neural networks used for the ratio estimation, and that in this simulation study the values $n^L_1=n^L_{2,1}=n^L_{2,2}=8$ and $n^N_1=n^N_{2,1}=n^LN_{2,2}=64$ used in Section \ref{subsec:SimulationStudy-FittingResults} provide good results. For the one-step estimation results, shallow but large networks might be slightly preferable to the chosen architecture, but not in a way sufficient to change our analysis. Indeed, it appears clearly that modifying the architecture and/or increasing the number of parameters in the one-step estimation is not sufficient to improve the estimation to the level of the two-step estimation, stressing the importance of taking advantage of the multiplying structure in the estimation.

\subsubsection{Parsimony and computational time}

We end this simulation study with a few comments on the parsimony and the computational cost of the estimation methods.
If we set $n^L_1=n^L_{2,1}=n^L_{2,2}$ and $n^N_1=n^N_{2,1}=n^LN_{2,2}$ as we did above, then the two-step estimation has a much larger number of parameters since we use $1+\#\mathbb I = 5$ networks very close in shape to the single one used for the one-step estimation method.
In the case $n^L_1=n^L_{2,1}=n^L_{2,2}=8$ and $n^N_1=n^N_{2,1}=n^LN_{2,2}=64$, this represents 33,991 parameters for the one-step estimation method and 167,559 parameters in the two-step estimation method. However, the networks of the two-step estimation are trained on smaller subsamples and convergence is attained quickly, so that in our example the total estimation time for the two-step estimation method is only approximately twice the time used for the one-step estimation. Moreover, since the multiple networks use in the two-step estimation can in fact be trained in parallel, the two-step estimation can in fact be faster than the one-step estimation.

\subsection{An application to high-frequency trades and LOB data}

In this section we use limit order book data of the stock Total Energies SA (ISIN : FR0000120171) traded in Euronext Paris. Our dataset covers 22 trading days, from January 2nd, 2017 to January 31st, 2017. For each trading day, the dataset lists all market and marketable orders (all referred to as market orders hereafter) submitted to the exchange between 9:05 and 17:25 local time, i.e. excluding a few minutes after the opening auction and before the closing auction. For each submission, the dataset lists the timestamp of the order with microsecond precision, as well the limit order book (LOB) data at the first level, namely best bid and ask prices and quantities. In the following, for an order entering the system at time $t$, $a(t-)$ (resp. $b(t-)$) is the ask price (resp. bid price) and $q^A(t-)$ (resp. $q^B(t-)$) is the quantity available at the best ask (resp. bid) queue of the limit order book just before $t$. If a market order triggered multiple transactions, then only one market order is in the dataset. The resulting number of market orders in the sample is greater than 1,750,000.

Let $N^0$ denote the counting process of market orders submitted on the bid side (sell market order) and $N^1$ denote the counting process of market orders submitted on the ask side (buy market order). Each order is marked $0$ if it does not change the mid-price, and marked $1$ if it changes the mid-price (which is equivalent to say that its execution depletes the best quote, or that the size of the order is greater that the size of the best quote). The {\ablue dataset} can thus be modeled by a point process $((N^{i,k_i}_t)_{t\geq 0})_{i=0,1,k=0,1}$ which can either be seen as a 4-dimensional point process or as a 2-dimensional point process with marks in $\{0,1\}$.

The Level-I order book data can be used to compute significant covariates in a high-frequency finance context. Let $X^0_{t-} := i(t-) := \frac{q^B(t-)-q^A(t-)}{q^B(t-)+q^A(t-)}$ the imbalance measured just before the submission of an order at time $t$. Imbalance is a well-known indicator of the short-term behaviour of the market: an imbalance close to 1 (resp. -1) indicate a positive (resp. negative) pressure on the price.
Let $X^1_{t-}$ be the sign of the last trade, i.e. $X^1_{t-}=-1$ if the last transaction occured on the bid side of the limit order book, $X^1_{t-}=1$ if the last transaction occured on the ask side of the limit order book. It is well-known in high-frequency finance that the series of trade signs have long-memory and are thus informative in our context. Finally, let $X^2_{t-}:=s(t-):=a(t-)-b(t-)$ be the bid-ask spread measured just before the submission of a market order at time $t$. When the spread is greater than 1 tick, a trader can gain priority by placing limit orders inside the bid-ask spread and thus get faster execution without using market orders. The spread is thus informative in an intensity model for the point process $((N^{i,k_i}_t)_{t\geq 0})_{i=0,1,k=0,1}$. In the following,  {\ablue the} spread is expressed in number of ticks and takes values in $\{1,2,3\}$ (in the rare cases (2\% of the dataset) where the spread is greater than 3 ticks, we set it equal to 3 ticks).

We can write two intensity models for the submission of market orders in a limit order book. The first intensity model is simply written
\begin{equation}
	\lambda^{i,k_i}(t) = \lambda_0(t)\, \lambda^{i,{\ablue k_i}}(X^0_t,X^1_t,X^2_t),
\end{equation}
with $i=0,1$ (bid side or ask side), $k_i=0,1$ (not price-changing or price-changing) and the c\`agl\`ad processes $X^j$, $j=1,2,3$ are defined above.
The second intensity model for the submission is written as the marked ratio model
\begin{equation}
	\lambda^{i,k_i}(t) = \lambda_0(t)\, \lambda^i(X^0_t,X^1_t,X^2_t)\,p_i^{k_i}(X^0_t,X^1_t,X^2_t).
\label{eq:marketIntensityWithMarks}
\end{equation}

The first model can be estimated with the one-step estimation method of Section \ref{subsec:oneStepEstimation}. The second model can be estimated with the two-step estimation methods of Section \ref{subsec:twoStepEstimation}. In both cases the neural networks are defined with parameters $n^L_1=n^L_{2,1}=n^L_{2,2}=8$ and $n^N_1=n^N_{2,1}=n^LN_{2,2}=64$.

\begin{figure}
\begin{tabular}{c}
\includegraphics[width=\textwidth]{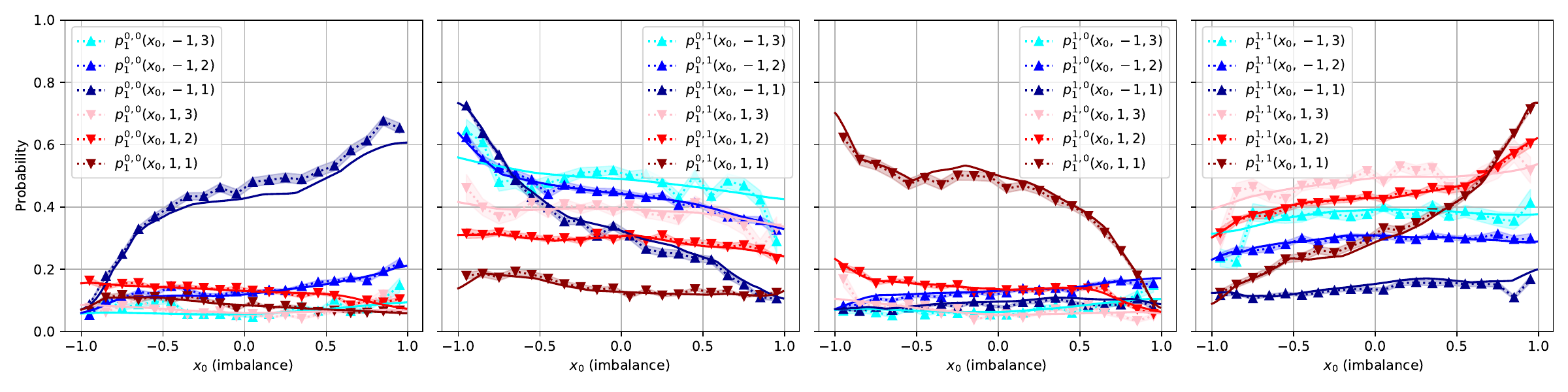}
\\
\includegraphics[width=\textwidth]{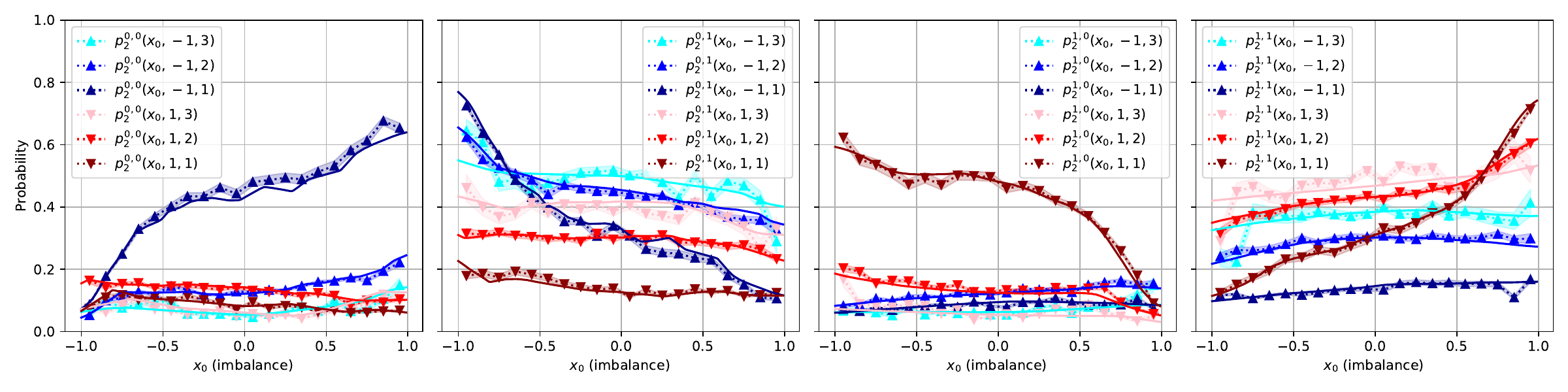}
\end{tabular}
\caption{Sign and price-changing character of trades -- Joint probabilities $p^{i,k}(x_0,x_1,x_2)$. One-step  (above) and two-step (below) estimation method. From left to right column-wise, $(i,k)=(0,0), (0,1), (1,0)$ and $(1,1)$. Empirical values with triangle and dotted lines. Fitted values as plain lines. Recall that $p^{i,k}(x_0,x_1,x_2)$ is the probability to observe a market order on the side $i$ with price-changing character $k$ when the LOB has imbalance $x_0$, spread $x_2$ and the last traded order was if sign $x_1$.)}
\label{fig:TradesSignAgg}
\end{figure}

Figure \ref{fig:TradesSignAgg} plots the fitting results. The first row plot the fitted probabilities $\hat p_1^{i,k_i}(x_0,x_1,x_2)$ (one-step estimation). From left to right, the four columns plot the probabilities in the case $(i,k_i)=((0,0)$, then $(0,1)$, then $(1,0)$ and finally $(1,1)$. On each plot, we have six curves corresponding to the cases $x_1,x_2\in\{-1,1\}\times \{1,2,3\}$. Imbalance $x_0$ is set as the abscissa of each plot. The second row provides the same plot for $\hat p_2^{i,k_i}(x_0,x_1,x_2)$, i.e. for the two-step estimation. 
It appears that both estimation methods give excellent fitting results. No method provides a strikingly better fit than the other. The model captures well-known {\ablue characteristics} of order flows in market microstructure: when the spread is equal to one tick, given that the previous order was a sell market order, the probability to observe another sell market order is high and
{\ablue the lesser the imbalance the higher the probability} 
that the order will deplete the best quote and change the price. When the spread increases, the curves flatten as the dependency to the imbalance is less strong. Observations are symmetric for buy market orders. A simple parametric form of the functionals $\lambda^{i,k}$, $ \lambda^i$, and $p_i^{k_i}$ would not be able to reproduce this variety of shapes (exponential forms have been tested and not shown here, because of poor results).
All in all, given these fitting results and the analysis of the simulation study, the multiplicative structure of Equation \eqref{eq:marketIntensityWithMarks} with a deep learning architecture seems well-suited for a intensity model for market orders depending in the observed spread, imbalance and last trade sign.

\section{{\bred Some basic estimates and proof of Theorem \ref{0607170208}}}\label{0610080715}
\subsection{Preparations}
First, we replace 
$R_T$ by its empirical version.  
Define the empirical error $\calr_T$ by
\beas 
\calr_T
&=& 
T^{-1}\bigg(-\int_0^T\big\{\wh{a}_T(X_t)-a^*(X_t)\big\}\cdot dN_t+\int_0^T\big\{\wh{b}_T(X_t)-b^*(X_t)\big\}dt\bigg), 
\eeas
and the expected empirical error of $(\wh{a}_T,\wh{b}_T)$ by 
\bea\label{0607160258}
R_T^e
&=& 
E\big[\calr_T]. 
\eea
The compensated $N$ is denoted by $\wt{N}$, that is, $d\wt{N}=dN-\lambda^*(X_t)dt$. 
For $T,T_1\in\bbT$, we have 
%
\bea\label{0607160122}
\big|R_T-R^e_T\big|
&\leq&
\big|\Phi^{(\ref{0512231520})}_T\big|+\big|\Phi^{(\ref{0512231521})}_T\big|, 
\eea
where
\bea\label{0512231520}
\Phi^{(\ref{0512231520})}_T
&=& 
T^{-1}E\bigg[-\int_0^T\big\{\wh{a}_T(X_t)-a^*(X_t)\big\}\cdot d\wt{N}_t\bigg]
\eea
and 
%
\bea\label{0512231521}
\Phi^{(\ref{0512231521})}_T
&=& 
{\colorr E\bigg[T_1^{-1}\int_0^{T_1}\wh{U}_T(\ol{X}_t)dt\bigg]}
\nn\\&&
-T^{-1}E\bigg[-\int_0^T\big\{\wh{a}_T(X_t)-a^*(X_t)\big\}\cdot \lambda^*(X_t)dt
+\int_0^T\big\{\wh{b}_T(X_t)-b^*(X_t)\big\}dt\bigg]
\nn\\&=&
T^{-1}E\bigg[\int_0^T\big\{\wh{U}_T(\ol{X}_t)-\wh{U}_T(X_t)\big\}dt\bigg]. 
\eea

For any $\delta>0$, {\bred we consider} a $\delta$-net 
$\big\{\{(a,b);\bbd((a,b),(a_k,b_k))<\delta\}\big\}_{k\in\calk_T}$ {\ared of ${\mathfrak F}_T$ }%
such that each ball has the radius $\delta$ in $\bbd$. 
{\bred We may assume that $\#\calk_T<\infty$; otherwise, the targeted inequality (\ref{202412220016}) is trivial. 
So we let $\calk_T=\{1,...,\caln_T\}$. As already mentioned, $\caln_T$ depends on $\delta$ as well as $T$.} %
We denote by $(a_{\tt k},b_{\tt k})$ 
the center of a $\delta$-ball for which the distance to $(\wh{a}_T,\wh{b}_T)$ is mimimum
among all the centers $(a_k,b_k){\bred\>\in\bbA}$ ($k=1,...,\caln_T$), 
where ${\tt k}$ is a random variable that indicates the number of one of the nearest points. 
The gap between {\ared $(a_{\tt k},b_{\tt k})$ and $(a^*,b^*)$} is evaluated with the function
\beas 
U^{\tt k}_T(x)
&=&
-\lambda^*(x)\cdot\big(a_{\tt k}(x)-a^*(x)\big)+\big\{b_{\tt k}(x)-b^*(x)\big\}. 
\eeas 
\begin{en-text}
For any $\delta>0$, there exists a $\delta$-net in the space of $(a,b)$ 
such that each ball has ``radius'' $\delta$. 
The ``radius'' is defined as $E[U(X_0)]$, or a distance equivalent to it. 
Denote by $\{(a_k,b_k)\}_k$ the centers of covering balls. 
\end{en-text}

\subsection{Estimate of $\Phi^{(\ref{0512231521})}_T$}
Define $r^k_T$ by 
\bea\label{0704012226}
r^k_T 
&=&
(\ttT^{-1}{\colorb(\log \ttT)^2}\log\caln_T)^{1/2}
\nn\\&&
\vee
\bigg({\tt h}^{-1}\int_0^{\tt h}E\bigg[-\lambda^*(\ol{X}_t)\cdot\big\{a_k(\ol{X}_t)-a^*(\ol{X}_t)\big\}
+\big\{b_k(\ol{X}_t)-b^*(\ol{X}_t)\big\}\bigg]dt\bigg)^{1/2}
\eea
for $k\in\{1,...,\caln_T\}$. 
The random number $r^{\tt k}_T$ is $r^k_T$ with ${\tt k}$ plugged into $k$. 

We have 
\bea\label{0512261023}
\big|\Phi^{(\ref{0512231521})}_T\big|
&=&
\bigg|T^{-1}E\bigg[\int_0^T\big\{\wh{U}_T(\ol{X}_t)-\wh{U}_T(X_t)\big\}dt\bigg]\bigg|
\nn\\&\leq&
\bigg|T^{-1}E\bigg[\int_0^T\big\{U^{\tt k}_T(\ol{X}_t)-U^{\tt k}_T(X_t)\big\}dt\bigg]\bigg|
+{\colorr\delta}.
\eea

The compatibility condition (\ref{202412111052}) implies
{\colorr
\bea\label{202412111732}&&
\bigg({\tt h}^{-1}\int_0^{\tt h}E_{\ol{X}}\bigg[-\lambda^*(\ol{X}_t)\cdot\big\{a_k(\ol{X}_t)-a^*(\ol{X}_t)\big\}
+\big\{b_k(\ol{X}_t)-b^*(\ol{X}_t)\big\}\bigg]dt\bigg)^{1/2}\bigg|_{k={\tt k}}
\nn\\&\leq&
{\bred\constzz}
\bigg({\bred\tth^{-1}}\int_0^{\tt h}E_{\ol{X}}\bigg[\big|a_k(\ol{X}_t)-a^*(\ol{X}_t)\big|^2
+\big|b_k(\ol{X}_t)-b^*(\ol{X}_t)\big|^2\bigg]dt\bigg)^{1/2}\bigg|_{k={\tt k}}
\nn\\&\leq&
{\bred\constzz}
\bigg({\bred\tth^{-1}}\int_0^{\tt h}E_{\ol{X}}\bigg[\big|\wh{a}_T(\ol{X}_t)-a^*(\ol{X}_t)\big|^2
+\big|\wh{b}_T(\ol{X}_t)-b^*(\ol{X}_t)\big|^2\bigg]dt\bigg)^{1/2}
\nn\\&&
+{\bred\constzz}\bigg({\bred\tth^{-1}}\int_0^{\tt h}E_{\ol{X}}\bigg[\big|\wh{a}_T(\ol{X}_t)-a_k(\ol{X}_t)\big|^2
+\big|\wh{b}_T(\ol{X}_t)-b_k(\ol{X}_t)\big|^2\bigg]dt\bigg)^{1/2}\bigg|_{k={\tt k}}
\nn\\&&\qquad(\text{by the triangular inequality})
\nn\\&\leq&
{\bred\constzz}\bigg({\bred\tth^{-1}}\int_0^{\tt h}E_{\ol{X}}\bigg[\big|\wh{a}_T(\ol{X}_t)-a^*(\ol{X}_t)\big|^2
+\big|\wh{b}_T(\ol{X}_t)-b^*(\ol{X}_t)\big|^2\bigg]dt\bigg)^{1/2}
+\delta
\nn\\&\leq&
{\bred\constzz^2}
\bigg({\tt h}^{-1}\int_0^{\tt h}E_{\ol{X}}\bigg[-\lambda^*(\ol{X}_t)\cdot\big\{\wh{a}_T(\ol{X}_t)-a^*(\ol{X}_t)\big\}
+\big\{\wh{b}_T(\ol{X}_t)-b^*(\ol{X}_t)\big\}\bigg]dt\bigg)^{1/2}
+\delta.
\nn\\&&
\eea
}
Then, 
\bea\label{0512241243}
r^{\tt k}_T
&\leq&(\ttT^{-1}{\colorb(\log \ttT)^2}\log\caln_T)^{1/2}\nn\\&&+\bigg({\tt h}^{-1}\int_0^{\tt h}E_{\ol{X}}\bigg[-\lambda^*(\ol{X}_t)\cdot\big\{a_k(\ol{X}_t)-a^*(\ol{X}_t)\big\}+\big\{b_k(\ol{X}_t)-b^*(\ol{X}_t)\big\}\bigg]dt\bigg)^{1/2}\bigg|_{k={\tt k}}
\nn\\&\leq^{(\ref{202412111732})}&
(\ttT^{-1}{\colorb(\log \ttT)^2}\log\caln_T)^{1/2}\nn\\&&
+{\bred\constzz^2}\bigg({\tt h}^{-1}\int_0^{\tt h}E_{\ol{X}}\bigg[-\lambda^*(\ol{X}_t)\cdot\big\{\wh{a}_T(\ol{X}_t)-a^*(\ol{X}_t)\big\}+\big\{\wh{b}_T(\ol{X}_t)-b^*(\ol{X}_t)\big\}\bigg]dt\bigg)^{1/2}+\delta 
\nn\\&=&
(\ttT^{-1}{\colorb(\log \ttT)^2}\log\caln_T)^{1/2}+{\bred\constzz^2}\wh{E}_T^{1/2}+\delta,
\eea
where 
\beas
\wh{E}_T
&=&
{\tt h}^{-1}\int_0^{\tt h}
E_{\ol{X}}
\big[{\bred\wh{U}_T(\ol{X}_t)}\big]
dt. 
\eeas
{\ared For simplicity of the presentation, we often write inequalities like $\leq^{(*.*)}$ indicating use of the item provided by $(*.*)$. }

\begin{en-text}
\im Try once again. ${\tt k}$ denotes the (random) number indicating the center of the ball 
nearest to ($\wh{a}_T,\wh{b}_T)$. 
\bea\label{0512261023}
\big|\Phi^{(\ref{0512231521})}_T\big|
&=&
\bigg|T^{-1}E\bigg[\int_0^T\big\{\wh{U}_T(\ol{X}_t)-\wh{U}_T(X_t)\big\}dt\bigg]\bigg|
\nn\\&\leq&
\bigg|T^{-1}E\bigg[\int_0^T\big\{U^{\tt k}_T(\ol{X}_t)-U^{\tt k}_T(X_t)\big\}dt\bigg]\bigg|
+{\colorr\delta}
\nn\\&\leq&
\bigg|E\bigg[\bigg(r_T^{\tt k}FT^{-1}\bigg)
\times (r_T^{\tt k}F)^{-1}\int_0^T\big\{U^{\tt k}_T(\ol{X}_t)-U^{\tt k}_T(X_t)\big\}dt\bigg]\bigg|
+{\colorr\delta}
\eea
where 
the factor $r_T^{\tt k}$ is a random tuning parameter determined by ${\tt k}$: 
there is a deterministic function $k\mapsto r^k_T$ and $r^{\tt k}_T=r^k_T|_{k={\tt k}}$. 
The explicit form will be specified below. 
The first factor on the right-hand side of the above inequality will be $R_T^{1/2}$. 
The second factor will be estimated by LD. 
\end{en-text}
\begin{en-text}
\im {\colorr $U^k(x)$
\beas 
U^k(x_{[0,{\tt h}]})
&=& 
{\tt h}^{-1}\int_0^{\tt h}\big[-\lambda^*(x_t)\cdot\big\{a_k(x_t)-a^*(x_t)\big\}+\big\{b_k(x_t)-b^*(x_t)\big\}\big]dt. 
\eeas
}
\end{en-text}

Let 
$
U^{\tt k}(x)=-\lambda^*(x)\cdot\big\{a_{\tt k}(x)-a^*(x)\big\}+\big\{b_{\tt k}(x)-b^*(x)\big\}
$ 
and 
\bea\label{0702130837}
\bbL_T
&=& 
(r_T^{\tt k})^{-1}F^{-1}{\ared\tth^{-1}}\int_0^T\big\{U^{\tt k}(\ol{X}_t)-U^{\tt k}(X_t)\big\}dt.
\eea
Then, by (\ref{0512261023}) and (\ref{0512241243}), 
\bea\label{0607160130}
\big|\Phi^{(\ref{0512231521})}_T\big|
&\leq&
\bigg|E\bigg[\big(r_T^{\tt k}F\ttT^{-1}\big)
\times (r_T^{\tt k}F{\ared\tth})^{-1}\int_0^T\big\{U^{\tt k}_T(\ol{X}_t)-U^{\tt k}_T(X_t)\big\}dt\bigg]\bigg|
+{\colorr\delta}
\nn\\&=&
\big|E\big[\big(r_T^{\tt k}F\ttT^{-1}\big)
\times\bbL_T\big]\big|
+{\colorr\delta}
\nn\\&\leq&
{\bred\constzz^2}F\ttT^{-1}{\ared\big|}E\big[\wh{E}_T^{1/2}\bbL_T\big]\big|
+F\ttT^{-1}{\ared\big|}E\big[\big\{(\ttT^{-1}{\colorb(\log \ttT)^2}\log\caln_T)^{1/2}+\delta\big\}\bbL_T\big]\big|
+{\colorr\delta}
\nn\\&\leq&
{\bred\constzz^2}F\ttT^{-1}R_T^{1/2}\big(E[|\bbL_T|^2]\big)^{1/2}
+F\ttT^{-1}\bigg\{(\ttT^{-1}{\colorb(\log \ttT)^2}\log\caln_T)^{1/2}+\delta\bigg\}E[|\bbL_T|]
+{\colorr\delta}.
\nn\\&&
\eea

{\colorr
\subsection{A large deviation estimate for an additive functional}\label{0704012250}
{\ared Recall that }%
the covariate process $X$ takes values in a measurable set $\calx$ in $\bbR^{\sfd_X}$. 
It is assumed that $X$ is periodically stationary. 
{\ared For a bounded measurable function $\ttU:\calx\to\bbR_+=[0,\infty)$, 
}
let 
\bea\label{202504131628}
\bbZ_\ell^{(T)}
&=& 
(r_T)^{-1}{\tt h}^{-1}
\int_{(\ell-1){\tt h}}^{\ell{\tt h}}{\bred\big\{\ttU(X_t)-E[\ttU(X_t)]\big\}}dt\qquad(\ell\in\bbN,\>T\in\bbT)
\eea
for 
\bea\label{0606221447}
r_T 
&=&
(\ttT^{-1}{\colorb(\log \ttT)^2}\log\caln_T)^{1/2}\vee\big(E\big[{\ared\ttU}(X_{[0,{\tt h}]})\big]\big)^{1/2},\qquad
{\ared\ttU}(X_{[0,{\tt h}]}) \yeq {\tt h}^{-1}\int_0^{\tt h}{\ared\ttU}(X_t)dt.
\eea
From (\ref{0606221447}), in particular, 
\bea\label{0606221449}
r_T\ygeq (\ttT^{-1}{\colorb(\log \ttT)^{\ared2}}\log\caln_T)^{1/2}, 
\quad\text{equivalently,}\quad
r_T^{-1}\yleq \ttT^{1/2}{\colorb(\log \ttT)^{-1}}(\log\caln_T)^{-1/2}
\eea
and 
\bea\label{0606221600}
r_T^2 &\geq& E\big[{\ared {\tt U}}(X_{[0,{\tt h}]})\big].
\eea

%
The following lemma gives a large deviation inequality for the sum $\sum_{\ell=1}^\ttT\bbZ_\ell^{(T)}$. 
\begin{lemma}\label{0606190218}
Let $\ep$ and $\consta$ be positive numbers. 
{\ared Suppose that }
\bea
\ttT&\geq& {\cred3}\vee\log\caln_T,
\label{0606221611}
\\
\log\caln_T&\geq&4\|{\ared {\tt U}}\|_\infty^2 \label{0606221457},
\\
x&\geq&\consta {\ared\ttT}^{1/2}(\log\caln_\ttT)^{1/2} \label{0606221439}. 
\eea
Then, for some positive constant $\constb$ depending only on ${\ared\gamma}$, it holds that 
\bea\label{0606190222}
P\bigg[\bigg|\sum_{\ell=1}^{\ttT}\bbZ_\ell^{(T)}\bigg|\geq x\bigg]
&\leq&
\exp\bigg[-\frac{\constb x^{\frac{1-\ep}{1+\ep}}}{K({\ared\gamma},\|{\ared\ttU}\|_\infty,\consta,\ep,T)}\bigg]\qquad(x>0,\>T\in\bbT)
\eea
where 
\bea\label{0606221729}&&
K({\ared\gamma},\|{\ared\ttU}\|_\infty,\consta,\ep,T)
\nn\\&=& 
(1+\|{\ared\ttU}\|_\infty)^2{\ablue(\consta^{-1}+1)}
\consta^{-\frac{2\ep}{1+\ep}}
\big(V({\ared\gamma},\ep)+{\colorb\log \ttT}\big)\ttT^{1/2}(\log\caln_T)^{-1/2-\frac{2\ep}{1+\ep}}
\eea
for a constant $V({\ared\gamma},\ep)$ given by 
\beas
V({\ared\gamma},\ep)
&=&
4\bigg[
1+4\sum_{j\in\bbN}{\ared\gamma}^{-\frac{1}{1+\ep^{-1}}}\exp\big(-\frac{{\ared\gamma}}{1+\ep^{-1}}j\big)\bigg]. 
\eeas
\end{lemma}
\proof 
{\ared We may assume that $\|{\ared\ttU}\|_\infty>0$; otherwise, 
the inequality (\ref{0606190222}) is trivial since $\bbZ_\ell^{(T)}=0$. 
Then $\caln_T>1$ by (\ref{0606221457}). }

From Theorem 2 of Merlev\`ede et al. \cite{merlevede2009bernstein}, we have 
\bea\label{060619024}
P\bigg[\bigg|\sum_{\ell=1}^{\ttT}\bbZ_\ell^{(T)}\bigg|\geq x\bigg]
&\leq&
\exp(-I_T(x))
\eea
for all $\ttT\geq{\cred3}$, where 
\bea\label{0606190241}
I_T(x) &=& \frac{\constb(r_T)^2x^2}{v^2\ttT+4\|{\ared\ttU}\|_\infty^2+2\|{\ared\ttU}\|_\infty(\log \ttT)^2r_Tx}.
\eea
The constant $\constb$ depends only on ${\ared\gamma}$, and the constant $v$ is given by 
\bea\label{0606190242}
v^2 
&=& 
\sup_{\ell\in\bbN}\bigg[\text{Var}[\wh{\bbZ}_\ell^{(T)}]+2\sum_{j>\ell}\big|\text{Cov}[\wh{\bbZ}_\ell^{(T)},\wh{\bbZ}_j^{(T)}]\big|\bigg]
\eea
for 
\beas 
\wh{\bbZ}_\ell^{(T)}
&=& 
{\tt h}^{-1}
\int_{(\ell-1){\tt h}}^{\ell{\tt h}}{\bred\big\{\ttU(X_t)-E[\ttU(X_t)]\big\}}dt. 
\eeas
%

The covariance inequality (Rio \cite{rio2017asymptotic} p. 6) gives 
\bea\label{0606190210}
v^2 
&\leq&
V({\ared\gamma},\ep)\big(E[{\ared\ttU}(X_{[0,\tth]}{\ared)}]\big)^{\frac{1}{1+\ep}}\>\|{\ared\ttU}\|_\infty^{\frac{1+2\ep}{1+\ep}}.
\eea
Remark that $\frac{1}{2+2\ep}+\frac{1}{2+2\ep}+\frac{1}{1+\ep^{-1}}=1$ and that 
\beas 
\|{\ared\ttU}(X_{[0,{\tt h}]})\|_{2+2\ep} 
&\leq& 
\big(E[{\ared\ttU}^{2+2\ep}(X_{[0,{\tt h}]})]\big)^{\frac{1}{2+2\ep}}
\yleq 
\big(E[{\ared\ttU}(X_{[0,{\tt h}]})]\big)^{\frac{1}{2+2\ep}}\>\|{\ared\ttU}(X_{[0,{\tt h}]})\|_\infty^{\frac{1+2\ep}{2+2\ep}}, 
\eeas
{\bred 
additionally, 
$\|E[{\ared\ttU}(X_{[0,{\tt h}]})]\|_{2+2\ep}= E[{\ared\ttU}(X_{[0,{\tt h}]})]\leq \big(E[{\ared\ttU}(X_{[0,{\tt h}]})]\big)^{\frac{1}{2+2\ep}}\>\|{\ared\ttU}(X_{[0,{\tt h}]})\|_\infty^{\frac{1+2\ep}{2+2\ep}}$. 
}
For the constant $V({\ared\gamma},\ep)$, we have
\bea\label{0606190210b}
V({\ared\gamma},\ep)
&\leq&
4+\frac{16{\ared\gamma}^{-\frac{1}{1+\ep^{-1}}}}{1-\exp\big(-\frac{{\ared\gamma}}{1+\ep^{-1}}\big)}. 
\eea

{\bred We know}
\bea
r_T^{-1}&\leq^{(\ref{0606221449})}& \ttT^{1/2}{\colorb(\log \ttT)^{-1}}(\log\caln_T)^{-1/2}\label{0606192026}
\yeq
\ttT^{1/2}(\log\caln_T)^{1/2}{\colorb(\log \ttT)^{-1}}(\log\caln_T)^{-1} 
\nn\\&\leq^{(\ref{0606221439})}& 
\consta^{-1}{\colorb(\log \ttT)^{-1}}(\log\caln_T)^{-1}x \label{0606192027}
\eea
and hence 
\bea\label{0606192034}
(r_T)^{-\frac{2\ep}{1+\ep}}
&\leq&
\consta^{-\frac{2\ep}{1+\ep}} x^{\frac{2\ep}{1+\ep}}{\colorb(\log \ttT)^{{\cred-}\frac{2\ep}{1+\ep}}}(\log\caln_T)^{-\frac{2\ep}{1+\ep}}
\nn\\&\leq^{{\ared(\ref{0606221611})}}&
\consta^{-\frac{2\ep}{1+\ep}} x^{\frac{2\ep}{1+\ep}}(\log\caln_T)^{-\frac{2\ep}{1+\ep}}.
\eea

We have 
\bea\label{0606221503}
(r_T)^2\ttT\>\geq^{(\ref{0606221449})}\>{\colorb(\log \ttT)^{{\ared2}}}\log\caln_T
\>\geq^{(\ref{0606221611})\atop(\ref{0606221457})}\>4\|{\ared\ttU}\|_\infty^2, 
\eea
\bea\label{0606221528}
(r_T)^2\ttT
&{\ared\leq^{(\ref{0606221611})}}&
{\ared(1+\|{\tt U}\|_\infty^2)}{\colorb(\log \ttT)}(r_T)^2\ttT^{1/2}\cdot \ttT^{1/2}
\nn\\&\leq^{(\ref{0606221439})}&
{\ared(1+\|{\tt U}\|_\infty^2)}{\colorb(\log \ttT)}(r_T)^2\ttT^{1/2}\consta^{-1}(\log\caln_T)^{-1/2}x
\eea
and
\bea\label{0606221508}
2\|{\ared\ttU}\|_\infty(\log \ttT)^2r_Tx
&=&
2\|{\ared\ttU}\|_\infty(\log \ttT)^2(r_T)^2(r_T)^{-1}x
\nn\\&\leq^{(\ref{0606221449})}&
2\|{\ared\ttU}\|_\infty(\log\ttT)(r_T)^2\ttT^{1/2}(\log\caln_T)^{-1/2}x. 
\eea
From (\ref{0606221503}), (\ref{0606221528}) and (\ref{0606221508}), we obtain 
\bea\label{06062211649}
4\|{\ared\ttU}\|_\infty^2+2\|{\ared\ttU}\|_\infty(\log \ttT)^2r_Tx
&\leq&
{\ared(1+\|{\tt U}\|_\infty)^2}{\colorb(\log \ttT)}(r_T)^2\ttT^{1/2}(\consta^{-1}+1)(\log\caln_T)^{-1/2}x. 
\eea
Since 
\beas
x(\log\caln_T)^{-1}
\>\geq^{(\ref{0606221439})}\>
\consta \ttT^{1/2}(\log\caln_T)^{-1/2}
\>\geq^{(\ref{0606221611})}\>
\consta, 
\eeas
we have 
\bea\label{0606221541}&&
4\|{\ared\ttU}\|_\infty^2+2\|{\ared\ttU}\|_\infty(\log \ttT)^2r_Tx
\nn\\&\leq^{(\ref{06062211649})}&
{\ared(1+\|{\tt U}\|_\infty)^2}{\colorb(\log \ttT)}(r_T)^2\ttT^{1/2}(\consta^{-1}+1)(\log\caln_T)^{-1/2}x
\cdot x^{\frac{2\ep}{1+\ep}}(\log\caln_T)^{-\frac{2\ep}{1+\ep}}\consta^{-\frac{2\ep}{1+\ep}}
\nn\\&=&
{\ared(1+\|{\tt U}\|_\infty)^2}(\consta^{-1}+1)\consta^{-\frac{2\ep}{1+\ep}}{\colorb(\log \ttT)}(r_T)^2\ttT^{1/2}(\log\caln_T)^{-1/2-\frac{2\ep}{1+\ep}}
x^{1+\frac{2\ep}{1+\ep}}.
\eea

Moreover, 
\bea\label{0606221615}
\ttT^{1/2}(\log\caln_T)^{-1/2-\frac{2\ep}{1+\ep}} x^{1+\frac{2\ep}{1+\ep}}
&=&
\ttT^{1/2}(\log\caln_T)^{-1/2}x\cdot x^{\frac{2\ep}{1+\ep}}(\log\caln_T)^{-\frac{2\ep}{1+\ep}}
\nn\\&\geq^{\ared(\ref{0606221439})}&
\consta \ttT\cdot x^{\frac{2\ep}{1+\ep}}(\log\caln_T)^{-\frac{2\ep}{1+\ep}}
\nn\\&\geq^{(\ref{0606192034})}&
\consta^{{\cred1+\frac{2\ep}{1+\ep}}} \ttT(r_T)^{-\frac{2\ep}{1+\ep}}.
\eea
Then
\bea\label{0606192050}
v^2\ttT
&\leq^{(\ref{0606190210})}&
V({\ared\gamma},\ep)\big(E[{\ared\ttU}(X_{[0,{\tt h}]})]\big)^{\frac{1}{1+\ep}}\>\|{\ared\ttU}\|_\infty^{\frac{1+2\ep}{1+\ep}}\ttT
\nn\\&\leq^{(\ref{0606221600})}&
V({\ared\gamma},\ep)(r_T)^2\|{\ared\ttU}\|_\infty^{\frac{1+2\ep}{1+\ep}}(r_T)^{\frac{-2\ep}{1+\ep}}\>\ttT
\nn\\&\leq^{(\ref{0606221615})}&
V({\ared\gamma},\ep)(1+\|{\ared\ttU}\|_\infty^2)
\consta^{-1{\cred-\frac{2\ep}{1+\ep}}}
(r_T)^2\ttT^{1/2}(\log\caln_T)^{-1/2-\frac{2\ep}{1+\ep}} x^{1+\frac{2\ep}{1+\ep}}. 
\eea

From (\ref{0606221541}) and (\ref{0606192050}), we obtain 
\bea\label{0606221633}&&
v^2{\ared\ttT}+4\|{\ared\ttU}\|_\infty^2+2\|{\ared\ttU}\|_\infty(\log {\ared\ttT})^2r_Tx
\nn\\&\leq&
V({\ared\gamma},\ep)(1+\|{\ared\ttU}\|_\infty^2)\consta^{-1{\cred-\frac{2\ep}{1+\ep}}}(r_T)^2\ttT^{1/2}(\log\caln_T)^{-1/2-\frac{2\ep}{1+\ep}} x^{1+\frac{2\ep}{1+\ep}}
\nn\\&&
+{\ared(1+\|{\tt U}\|_\infty)^2}(\consta^{-1}+1)\consta^{-\frac{2\ep}{1+\ep}}{\colorb(\log {\ared\ttT})}(r_T)^2{\ared\ttT}^{1/2}(\log\caln_T)^{-1/2-\frac{2\ep}{1+\ep}}
x^{1+\frac{2\ep}{1+\ep}}
\nn\\&\leq&
(1+\|{\ared\ttU}\|_\infty)^2{\ablue(\consta^{-1}+1)}
\consta^{-\frac{2\ep}{1+\ep}}
\big(V({\ared\gamma},\ep)+{\colorb\log {\ared\ttT}}\big)
(r_T)^2{\ared\ttT}^{1/2}(\log\caln_T)^{-1/2-\frac{2\ep}{1+\ep}}
x^{1+\frac{2\ep}{1+\ep}}
\nn\\&=&
K({\ared\gamma},\|{\ared\ttU}\|_\infty,\consta,\ep,T)(r_T)^2x^{1+\frac{2\ep}{1+\ep}}. 
\eea

\begin{en-text}
\koko we have 
\bea\label{0606190343}&&
v^2T+4\|U\|_\infty^2+2\|U\|_\infty(\log T)^2r_Tx
\nn\\&\leq&
V({\ared\gamma},\ep)\big(E[U(X_0)]\big)^{\frac{1}{1+\ep}}\>\|U\|_\infty^{\frac{1+2\ep}{1+\ep}}T
+4\|U\|_\infty^2+2\|U\|_\infty(\log T)^2r_Tx
\nn\\&\leq&
V({\ared\gamma},\ep)(r_T)^{\frac{2}{1+\ep}}\>\|U\|_\infty^{\frac{1+2\ep}{1+\ep}}T+(r_T)^2T+2\|U\|_\infty(\log T)^2r_Tx
\nn\\&&\hspace{150pt}(\koko\because (r_T)^2T\ygeq\log\caln_T\ygeq^{(\ref{0606221457})}4\|U\|_\infty^2)
\nn\\&=&
V({\ared\gamma},\ep)(r_T)^2(r_T)^{\frac{-2\ep}{1+\ep}}\>\|U\|_\infty^{\frac{1+2\ep}{1+\ep}}T+(r_T)^2T+2\|U\|_\infty(\log T)^2(r_T)^2(r_T)^{-1}x
\nn\\&\leq^{(\ref{0606192026})}&
V({\ared\gamma},\ep)(r_T)^2(r_T)^{\frac{-2\ep}{1+\ep}}\>\|U\|_\infty^{\frac{1+2\ep}{1+\ep}}T+(r_T)^2T+2\|U\|_\infty(\log T)^2(r_T)^2T^{1/2}(\log\caln_T)^{-1/2}x\koko
\nn\\&\leq&
V({\ared\gamma},\ep)(r_T)^2(r_T)^{\frac{-2\ep}{1+\ep}}\>\|U\|_\infty^{\frac{1+2\ep}{1+\ep}}T+4\|U\|_\infty(\log T)^2(r_T)^2T^{1/2}(\log\caln_T)^{-1/2}x
\nn\\&& 
((r_T)^2T\leq 2\|U\|_\infty(\log T)^2(r_T)^2T^{1/2}(\log\caln_T)^{-1/2}x\quad\text{when}\quad 2\|U\|_\infty(\log T)^2\geq1)\koko
\nn\\&\leq^{(\ref{0606192034})}&
V({\ared\gamma},\ep)(r_T)^2(r_T)^{\frac{-2\ep}{1+\ep}}\>\|U\|_\infty^{\frac{1+2\ep}{1+\ep}}T
\nn\\&&
+4\|U\|_\infty(\log T)^2(r_T)^2T^{1/2}(\log\caln_T)^{-1/2}x\cdot x^{\frac{2\ep}{1+\ep}}(\log\caln_T)^{-\frac{2\ep}{1+\ep}}
\nn\\&&(x^{\frac{2\ep}{1+\ep}}(\log\caln_T)^{-\frac{2\ep}{1+\ep}}\geq1\quad\text{since}\quad
x\geq T^{1/2}(\log\caln_T)^{1/2}\quad\text{and}\quad \log\caln_T\leq T^{1-\ep_1})
\nn\\&\leq^{(\ref{0606192050})}&
2(1+\|U\|_\infty)^2(\log T)^2(r_T)^2T^{1/2}(\log\caln_T)^{-1/2}x\cdot x^{\frac{2\ep}{1+\ep}}(\log\caln_T)^{-\frac{2\ep}{1+\ep}}
\label{0606192052}
\nn\\&=&
2(1+\|U\|_\infty)^2(\log T)^2(r_T)^2T^{1/2}(\log\caln_T)^{-\half-\frac{2\ep}{1+\ep}}x^{\frac{1+3\ep}{1+\ep}}
\eea
\end{en-text}
\begin{en-text}
\bea\label{0606192031}
(r_T)^{-1}
&\leq^{(\ref{0606192026})}& 
T^{1/2}(\log\caln_T)^{-1/2}
\eea
\end{en-text}
Now, from (\ref{0606190241}) and (\ref{0606221633}), we obtain 
\bea\label{0606220343}
I_T(x) 
&\geq& 
\frac{\constb(r_T)^2x^2}{K({\ared\gamma},\|{\ared\ttU}\|_\infty,\consta,\ep,T)(r_T)^2x^{1+\frac{2\ep}{1+\ep}}}
\yeq
\frac{\constb x^{\frac{1-\ep}{1+\ep}}}{K({\ared\gamma},\|{\ared\ttU}\|_\infty,\consta,\ep,T)}. 
\eea
This completes the proof. 
\qed\halflineskip

\begin{lemma}\label{0606270208}
Let $\constg$, $\consth$, $\consti$ and $x$ be positive numbers {\ared with $\constg,{\bred\consti}\geq1$}. 
Suppose that 
\bea&&
\ttT\ygeq {\cred3}\vee\log\caln_T,
\label{z0606221611}
\\&&
\log\caln_T\ygeq4
\|{\ared\ttU}\|_\infty^2 \label{z0606221457}
\\&&
\constg \ttT^{\ared\consth/2} \ygeq x\ygeq\consti{\colorb(\log \ttT)} \ttT^{1/2}(\log\caln_T)^{1/2}. \label{z0606221439}
\eea
Then, for some positive constant $\constk$ 
{\bred depending on ${\ared\gamma}$}, 
{\bred it holds that}
\bea\label{z0606190222}
P\bigg[\bigg|\sum_{\ell=1}^{T/{\tt h}}\bbZ_\ell^{(T)}\bigg|\geq x\bigg]
&\leq&
\exp\bigg[-\frac{{\colorb(\constg)^{-1}e^{{\cred-}\consth}}\constk x}{{\ared K_0(\|{\ared\ttU}\|_\infty,T)}}\bigg],
\eea
where 
\bea\label{z0606221729}
{\ared K_0(\|{\ared\ttU}\|_\infty,T)}
&=& 
(1+\|{\ared\ttU}\|_\infty)^2
{\colorb(\log \ttT)}\ttT^{1/2}(\log\caln_T)^{-1/2}.
\eea
\end{lemma}
\proof
{\bred We may assume that $\|\ttU\|_\infty>0$. }%
Let 
$\ep=1/\log \ttT$ and $z=\consti{\colorb\log \ttT}$. 
Then (\ref{0606190210b}) gives the estimate
\bea\label{0606270325}
V({\ared\gamma},\ep) &\leq& \constj\log \ttT
\eea
for some constant $\constj$ only depending on ${\ared\gamma}$, and 
it follows from (\ref{0606221729}) and (\ref{z0606221611}) that 
\bea\label{0606270342}
K({\ared\gamma},\|{\ared\ttU}\|_\infty,\consta,\ep,T)
&\leq&
\constl (1+\|{\ared\ttU}\|_\infty)^2
{\colorb(\log \ttT)}\ttT^{1/2}(\log\caln_T)^{-1/2}
\eea
for some constant $\constl$ depending on {\bred$\gamma$}. 
We remark that 
{\cred
\beas 
(\log\caln_T)^{-\frac{2\ep}{1+\ep}}
\yleq
\big(1/\log2)^{\frac{2\ep}{1+\ep}}
\yleq 
\big(1/\log2)^{\frac{2}{\log3+1}}
\eeas
since $\caln_T\geq2$ from $\log\caln_T>0$. 
}
Moreover, 
\beas 
x^{\frac{-2\ep}{1+\ep}}
&\geq&
\big(\constg {{\ared\tt T}}^{{\ared\consth/2}}\big)^{\frac{-2\ep}{1+\ep}}
\ygeq 
{\ared\exp\bigg[-\frac{2\log\constg}{1+\log3}-\consth\bigg]}
\nn\\&\geq&
\exp\big[-(\log\constg+\consth)\big]\yeq \constg^{-1}e^{-\consth}. 
\eeas
Now Lemma \ref{0606190218} provides the inequality (\ref{z0606190222}). 
\qed\halflineskip
%
%

\subsection{{\bred Estimation of $E[|\bbL_T|^2]$}}
{\bred We write $\|\lambda^*\|_\infty$ for $\big\||\lambda^*|\big\|_\infty$. }%
Define $\bbL_T^{k}$ by 
\beas 
\bbL_T^{k}
&=&
(r_T^{k})^{-1}{\ared F^{-1}}{\tt h}^{-1}\int_0^T\big\{U^{k}(\ol{X}_t)-U^{k}(X_t)\big\}dt.
\eeas

\begin{lemma}\label{0703241528}
Suppose that 
{\ablue
{\bred$\ttT\geq3\vee\log\caln_T$ and $\caln_T\geq2$}.
}%
Then, there exists a constant $\consto$ such that 
\beas 
E\big[|\bbL_T|^2\big]
\yleq
E\big[(\max_k|\bbL_T^k|)^2\big]
\yleq
\consto\> {\bred\ttT(\log \ttT)^2}\log\caln_T
\eeas
for all $T\in\bbT$. 
\end{lemma}
\proof
{\bred 
Since $-\bbL_T^{k}$ is the sum of the integrals 
$\wt{\bbL}_T^k:=(r_T^{k})^{-1}{\ared F^{-1}}{\tt h}^{-1}\int_0^T\big\{U^{k}(X_t)-E[U^{k}(X_t)]\big\}dt$ and 
$-(r_T^{k})^{-1}{\ared F^{-1}}{\tt h}^{-1}\int_0^T\big\{U^{k}(\ol{X}_t)-E[U^{k}(\ol{X}_t)]\big\}dt$, 
it is sufficient to estimate $E[(\max_k|\wt{\bbL}_T^k|)^2]$ only. 
}
{\bred Let $\zeta=2^{-1}(\log2)^{1/2}$. 
Set $\ttU= {\bred\zeta\big(2F(\|\lambda^*\|_\infty+1)\big)^{-1}}U^k$ and $r_T=r_T^k$, then}
\beas 
\bigg|\sum_{\ell=1}^{T/{\tt h}}\bbZ^{(T)}_\ell\bigg|
&\leq& 
2{\ared\ttT^{3/2}}
{\colorb(\log\ttT)}^{{\ared-1}}(\log\caln_T)^{-1/2}\|{\ared\ttU}\|_\infty
\yleq
{\ared\ttT^{3/2}}
\eeas
{\bred for $\bbZ^{(T)}_\ell$ of (\ref{202504131628}) 
since $\|\ttU\|_\infty\leq\zeta$, $\ttT\geq3$ and $\caln_T\geq2$. }%
{\bred Let $\constn$ be an arbitrary positive number. }%
Then
\beas 
E\big[(\max_k|{\bred\wt{\bbL}_T^k}|)^2\big]
&=&
\int_0^{{\bred2(\|\lambda^*\|_\infty+1)\ttT^{3/2}}}2xP\big[\max_k{\bred\wt{\bbL}_T^k}>x\big]dx
\nn\\&\leq&
(\constn)^2 {\bred\ttT(\log \ttT)^2}\log\caln_T
\nn\\&&
+\caln_T\int_{\constn {\bred(\log \ttT)\ttT^{1/2}}(\log\caln_T)^{1/2}}^{{\bred\infty}}2x
\exp\bigg[-\frac{{\colorb(\constg)^{-1}e^{{\cred-}\consth}}\constk x{\bred\big(2(\|\lambda^*\|_\infty+1)\big)^{-1}}}
{{\ared K_0(\|{\ared\ttU}\|_\infty,T)}}\bigg]dx
\eeas
{\bred by Lemma \ref{0606270208}.}

{\bred Let $\constkp=\constk\big(2(\|\lambda^*\|_\infty+1)\big)^{-1}$. }%
Using Lemma \ref{0606231547} below, we obtain
\beas &&
\caln_T\int_{\constn {\bred(\log \ttT)\ttT^{1/2}}(\log\caln_T)^{1/2}}^{{\bred\infty}}2x
\exp\bigg[-\frac{{\colorb(\constg)^{-1}e^{{\cred-}\consth}}{\bred\constkp} x}{{\ared K_0(\|{\ared\ttU}\|_\infty,T)}}\bigg]dx
\nn\\&\leq&
{\cred\frac{\constv(1+\|{\ared\ttU}\|_\infty)^4}{(\constg)^{-2}e^{{\cred-2}\consth}{\bred\constkp}^2 }\>(\log {\bred\ttT})^2T(\log\caln_T)^{-1}}
\caln_T
\exp\bigg[-\frac{{\colorb(\constg)^{-1}e^{{\cred-}\consth}}{\bred\constkp} \constn \log\caln_T
}{2(1+\|{\ared\ttU}\|_\infty)^2}\bigg]
\nn\\&\leq&
{\bred\constw (\log \ttT)^2\ttT}
\eeas
{\cred with some constants $\constv$ and $\constw$}, 
suppose that $\constn$ is chosen so large as
\beas 
\frac{{\colorb(\constg)^{-1}e^{{\cred-}\consth}}{\bred\constkp} \constn}{{\bred8}
} &>& 1.
\eeas
This completes the proof. 
\qed\halflineskip

\subsection{Estimate of $\Phi^{(\ref{0512231520})}_T$}
We will estimate $\Phi^{(\ref{0512231520})}_T$. 
{\bred The constant $r_T^k$ is defined by (\ref{0704012226}). }
%
%
Since
\beas
\bigg|T^{-1}E\bigg[\int_0^T\big\{\wh{a}_T(X_t)-a_{\tt k}(X_t)\big\}\cdot d\wt{N}_t\bigg]\bigg|
&\leq&
E\big[T^{-1}{\ared\sum_{i\in\bbI}(N_T^i+\|(\lambda^i)^*\|_\infty T)}\big]\big\|\wh{a}_T-a_{\tt k}\big\|_\infty 
\nn\\&\leq&
2{\ared\sfd_N}\|\lambda^*\|_\infty\big\|\wh{a}_T-a_{\tt k}\big\|_\infty 
\nn\\&\leq&
\bbd\big((\wh{a}_T,\wh{b}_T),(a_{\tt k},b_{\tt k})\big)
\nn\\&\leq&
\delta, 
\eeas
we obtain 
{\bred 
\bea\label{0703240439}
\big|\Phi^{(\ref{0512231520})}_T\big|
\begin{en-text}
&=&
\bigg|T^{-1}E\bigg[-\int_0^T\big\{\wh{a}_T(X_t)-a^*(X_t)\big\}\cdot d\wt{N}_t\bigg]\bigg|
\nn\\&\leq&
\bigg|T^{-1}E\bigg[\int_0^T\big\{a_T^{\tt k}(X_t)-a^*(X_t)\big\}\cdot d\wt{N}_t\bigg]\bigg|
\nn\\&&
+\bigg|T^{-1}E\bigg[\int_0^T\big\{\wh{a}_T(X_t)-a^{\tt k}(X_t)\big\}\cdot d\wt{N}_t\bigg]\bigg|
\end{en-text}
&\leq&
\bigg|T^{-1}E\bigg[\int_0^T\big\{a_{\tt k}(X_t)-a^*(X_t)\big\}\cdot d\wt{N}_t\bigg]\bigg|
+{\colorr\delta}
\begin{en-text}
\nn\\&\leq&
E\bigg[\bigg|T^{-1}\int_0^T\big\{a_T^k(X_t)-a^*(X_t)\big\}\cdot d\wt{N}_t\big|_{k={\tt k}}\bigg|\bigg]
+{\colorr\delta}
\nn\\&\leq&
\sum_kE\bigg[\bigg|T^{-1}\int_0^T\big\{a_T^k(X_t)-a^*(X_t)\big\}\cdot d\wt{N}_t\bigg|\bigg]
+{\colorr\delta}
\end{en-text}
\nn\\&=&
\bigg|E\bigg[\big(r_T^{\tt k}F\ttT^{-1}\big)\times\big(r_T^{\tt k}F\tth\big)^{-1}\int_0^T\big\{a_{\tt k}(X_t)-a^*(X_t)\big\}\cdot d\wt{N}_t\bigg]\bigg|+\delta
\nn\\&\leq^{(\ref{0512241243})}&
E\bigg[
F\ttT^{-1}\big\{
\big(\ttT^{-1}(\log\ttT)^2\log\caln_T\big)^{1/2}+{\bred\constzz^2}\wh{E}_T^{1/2}+\delta\big\}|\bbM_T|\bigg]+\delta
\nn\\&\leq&
{\bred\constzz^2}F\ttT^{-1}R_T^{1/2}\big(E[|\bbM_T|^2]\big)^{1/2}
+F\ttT^{-1}\bigg\{\big(\ttT^{-1}(\log\ttT)^2\log\caln_T\big)^{1/2}+\delta\bigg\}E[|\bbM_T|]+\delta,
\nn\\&&
\eea
}
where 
\beas
\bbM_T 
&=& 
(r_T^{\tt k})^{-1}{\bred F^{-1}\tth^{-1}}\int_0^T\big\{a_{\tt k}(X_t)-a^*(X_t)\big\}\cdot d\wt{N}_t. 
\eeas

\subsection{{\bred Estimation of $E[|\bbM_T|^2]$}}
Let 
\beas
\bbM_T^k
&=& 
(r_T^k)^{-1}{\bred F^{-1}\tth^{-1}}\int_0^T\big\{a_k(X_t)-a^*(X_t)\big\}\cdot d\wt{N}_t.
\eeas
%
{\ared 
The terminal value of the predictable quadratic variation of the local martingale associated with $\bbM_T^k$ is 
\beas 
\calv_k(T)
&=&
{\bred (r_T^k)^{-2}F^{-2}\tth^{-2}}
\bigg\langle
\int_0^\cdot\big\{a_k(X_t)-a^*(X_t)\big\}\cdot d\wt{N}_t\bigg\rangle_T
\nn\\&=&
{\bred (r_T^k)^{-2}F^{-2}\tth^{-2}}
\int_0^T\sum_{i\in\bbI}\big\{(a_k(X_t)-a^*(X_t))^i\big\}^2\big(\lambda^*(X_t)\big)^idt
\eeas
since there are no common jumps. 
Then 
\beas 
\calv_k(T)
&\leq&
{\bred (r_T^k)^{-2}F^{-2}\tth^{-2}}
\sfd_N\|\lambda^*\|_\infty\int_0^T|a_k(X_t)-a^*(X_t)|^2dt
\nn\\&\simleq&
{\bred (r_T^k)^{-2}F^{-2}\tth^{-2}}
\int_0^T\bigg[-\lambda^*(X_t)\cdot\{a_k(X_t)-a^*(X_t)\}+\{b_k(X_t)-b^*(X_t)\}\bigg]dt
\eeas
due to the compatibility (\ref{0703240806}).  
Therefore, 
\bea\label{0704020039}
E\big[\calv_k(T)\big]
&\leq&
2^{-1}\constr{\bred F^{-2}\tth^{-1}\ttT},
\eea
where $\constr$ is a positive constant depending on 
$\sfd_N$, {\bred$\|\lambda^*\|_\infty$ and $\constzz$, }%
and independent of $k$. 
{\bred It is possible to enlarge $\constr$ as we like. 

Let 
\beas 
\wt{\calv}_k(T) &=& r_T^k\tth (\|\lambda^*\|_\infty)^{-1}\big(\calv_k(T)-E[\calv_k(T)]\big). 
\eeas
Then, by using $\bbZ^{(T)}_\ell$ of Section \ref{0704012250}, the functional $\wt{\calv}_k(T)$ can be represented as
\beas 
\wt{\calv}_k(T) &=&4\zeta^{-1}\sum_{\ell=1}^\ttT \bbZ^{(T)}_\ell
\eeas
with $r_T=r_T^k$ and 
$\ttU(x)=4^{-1}F^{-2}(\|\lambda^*\|_\infty)^{-1}\zeta\sum_{i\in\bbI}\big\{(a_k(x)-a^*(x))^i\big\}^2\big(\lambda^*(x)\big)^i$. 

We suppose that 
$\ttT\geq3\vee\log\caln_T$ and $\caln_T\geq2$. 
Let $x_T = 
2^{-1}\constr (\|\lambda^*\|_\infty)^{-1}F^{-2}r_T^k\ttT
$, then 
\beas 
x_T &\geq&
2^{-1}\constr (\|\lambda^*\|_\infty)^{-1}F^{-2}(\log \ttT)(\log\caln_T)^{1/2}\ttT^{1/2}. 
\eeas
{\bred 
Choose the positive numbers $\constg$ and $\consth$ (after setting $\constr$), so as 
$\constg\ttT^{\consth/2}\geq x_T\log\ttT$ whenever 
$\ttT\geq3\vee\log\caln_T$ and $\caln_T\geq2$. 
}%
%
Moreover, take a sufficiently large $\constr$ such that 
\bea\label{0704012348}
2^{-1}\constr(\|\lambda^*\|_\infty)^{-1}F^{-2} &\geq& 1=:\consti.
\eea
Let $\Omega(k,T)=\big\{\calv_k(T)\leq\constr F^{-2}\tth^{-1} (\log\ttT)\ttT\big\}$ for $z\geq1$. 
Then Lemma \ref{0606270208} gives
\bea\label{0704012343}
P\big[\Omega(k,T)^c\big]
&\leq^{(\ref{0704020039})}&
P\bigg[\calv_k(T)-E[\calv_k(T)]\geq 2^{-1}\constr F^{-2}\tth^{-1} (\log\ttT)\ttT\bigg]
\nn\\&\leq&
P\bigg[\big|\wt{\calv}_k(T)\big|\geq (\log\ttT)x_T\bigg]
\yeq 
P\bigg[\bigg|\sum_{\ell=1}^\ttT\bbZ^{(T)}_\ell\bigg|\geq4^{-1} \zeta(\log\ttT)x_T\bigg]
\nn\\&\leq&
\exp\bigg[-\frac{\constab\constr F^{-2}(\|\lambda^*\|_\infty)^{-1} (\log\ttT)
\log\caln_T}
{16}\bigg]
\eea
for $\constab=2^{-1}{\colorb(\constg)^{-1}e^{{\cred-}\consth}}\constk\zeta$ depending on $\gamma$. 

Due to e.g. Inequality 1 of Shorack and Wellner \cite{shorack2009empirical}, p.899, we obtain
\bea\label{0704112009}
P\big[\big|\bbM_T^k\big|\geq x,\>\Omega(k,T)
\big]
&\leq&
2\exp\bigg[-\frac{x^2}{2\constr  F^{-2}\tth^{-1} (\log\ttT)\ttT}
\psi\bigg(\frac{2F^2x}{\constr r_T^k (\log\ttT)\ttT}\bigg)\bigg]\quad(x>0)
\nn\\&&
\eea
for any $k$, where 
$\psi(y)=2y^{-2}\big[(y+1)\{\log(y+1)-1\}+1\big]$. 
%
\begin{en-text}
When $x\geq\constn (\log T)T^{1/2}(\log\caln_T)^{1/2}$, 
\beas 
\frac{r^k_Tx}{T}
&\geq&
T^{-1}(T^{-1}{\colorb(\log T)}\log\caln_T)^{1/2}\times (\log T)T^{1/2}(\log\caln_T)^{1/2}
\nn\\&=&
T^{-1}(\log T)^{3/2}\log\caln_T
\eeas
\end{en-text}

The {\ared second-order} moment of $\bbM_T=\bbM_T^{\tt k}$ 
is estimated as follows: 
{\bred 
\bea\label{0610061321}
E\big[\big|\bbM_T\big|^2\big]
&=&
\int_0^\infty 2xP\big[\big|\bbM_T\big|\geq x\big]dx
\nn\\&\leq&
\big[\constn \ttT^{1/2}(\log \ttT)(\log\caln_T)^{1/2}\big]^2
\nn\\&&
+\caln_T\sup_k\int_{\constn \ttT^{1/2}(\log \ttT)(\log\caln_T)^{1/2}}^\infty 2xP\big[\big|\bbM_T^k\big|\geq x\big]dx
\nn\\&\leq&
\big[\constn \ttT^{1/2}(\log \ttT)(\log\caln_T)^{1/2}\big]^2
\nn\\&&
+\caln_T\sup_k\int_{\constn \ttT^{1/2}(\log \ttT)(\log\caln_T)^{1/2}}^\infty 2xP\big[\big|\bbM_T^k\big|\geq x,\>\Omega(k,T)\big]dx
\nn\\&&
+\caln_T\sup_kE\big[\big|\bbM_T^k\big|^21_{\Omega(k,T)^c}\big]. 
\eea
}

Apply the Burkholder-Davis-Gundy inequality to obtain
\beas 
E\big[|\bbM^{k}_T|^4\big]
&\simleq&
(r_T^k)^{-4}{\bred F^{-4}\tth^{-4}}\sum_{i\in\bbI}E\bigg[\bigg(\int_0^T\big|\big(a_k(X_t)-a^*(X_t)\big)^i\big|^2dN_t^i\bigg)^2\bigg]
\nn\\&\leq&
2(r_T^k)^{-4}{\bred F^{-4}\tth^{-4}}\bigg\{
\sum_{i\in\bbI}E\bigg[\bigg(\int_0^T\big|\big(a_k(X_t)-a^*(X_t)\big)^i\big|^2d\wt{N}_t^i\bigg)^2\bigg]
\nn\\&&
+
\sum_{i\in\bbI}E\bigg[\bigg(\int_0^T\big|\big(a_k(X_t)-a^*(X_t)\big)^i\big|^2(\lambda^*(X_t))^idt\bigg)^2\bigg]
\bigg\}
\nn\\&\leq&
32(r_T^k)^{-4}\tth^{-4}\sfd_N(1+\|\lambda^*\|_\infty)^2(T+T^2)
\nn\\&\leq&
64(r_T^k)^{-4}\tth^{-2}(1+\tth^{-1})\sfd_N(1+\|\lambda^*\|_\infty)^2\ttT^2
\nn\\&\leq&
64\tth^{-2}(1+\tth^{-1})\sfd_N(1+\|\lambda^*\|_\infty)^2\ttT^4. 
\eeas
We then have 
\bea\label{0703240749}
E\big[|\bbM^{k}_T|^21_{\Omega(k,T,z)^c}\big]
&\leq&
E\big[|\bbM^{k}_T|^4\big]^{1/2}P\big[\Omega(k,T,z)^c\big]^{1/2}
\nn\\&\simleq^{(\ref{0704012343})}&
8\tth^{-1}(1+\tth^{-1})^{1/2}\sfd_N(1+\|\lambda^*\|_\infty)\ttT^2
\nn\\&&\times
\exp\bigg[-\frac{\constab\constr F^{-2}(\|\lambda^*\|_\infty)^{-1}(\log\ttT)
\log\caln_T}
{32}\bigg]
\eea

We set $\constr=\constad F^2$, and choose a sufficiently large $\constad$ so that 
\beas 
\constad
&\geq& 
\max\{2,32\constab^{-1}\}\|\lambda^*\|_\infty,
\eeas
additionally to (\ref{0704020039}). 
Then we obtain 
\bea\label{0704042130}
\caln_T\sup_kE\big[\big|\bbM_T^k\big|^21_{\Omega(k,T,z)^c}\big]
&\leq&
\constac
\eea
for some constant $\constac$ depending on $\gamma$, $\tth$, $\sfd_N$, $\constzz$ and $\|\lambda^*\|_\infty$, 
by using 
\beas 
\frac{\log\ttT}{\log3}\frac{\log\caln_T}{\log2}\geq\frac{\log\ttT}{\log3}+\frac{\log\caln_T}{\log2}-1
\eeas
due to $\ttT\geq3$ and $\caln_T\geq2$.

We will show
\bea\label{202412141033}
\frac{x^2}{2\constr F^{-2}\tth^{-1}(\log\ttT)\ttT}
\psi\bigg(\frac{2F^2x}{\constr (\log\ttT)r_T^k\ttT}\bigg)
&\geq&
\consts \ttT^{-1/2}(\log\caln_T)^{1/2}x
\eea
for all $x$ satisfying 
\bea\label{202412141017}
x&\geq&\constn (\log \ttT)\ttT^{1/2}(\log\caln_T)^{1/2},
\eea
where $\consts$ is a constant depending on $\tth$, $\sfd_N$, $\constzz$ and $\|\lambda^*\|_\infty$, but 
independent of $\constn\geq1$. 
 {\cred
First, we see 
\bea\label{202412132308}
\constu &:=& \inf_{y>0}(y+1)\psi(y) \>>\>0. 
\eea
}%
When $\frac{2F^2x}{\constr (\log \ttT)r_T^k\ttT}\leq1$, 
\beas
\frac{x^2}{2\constr (\log \ttT)F^{-2}\tth^{-1}\ttT}
\psi\bigg(\frac{2F^2x}{\constr (\log \ttT)r_T^k\ttT}\bigg)
&\geq^{(\ref{202412132308})}&
\frac{\constu x^2}{4\constad\tth^{-1}\ttT\log\ttT}
\nn\\&\geq^{(\ref{202412141017})}&
\frac{\constu\constn  }{4\constad\tth^{-1}}\ttT^{-1/2}(\log\caln_T)^{1/2}x
\nn\\&\geq&
\frac{\constu }{4\constad\tth^{-1}}\ttT^{-1/2}(\log\caln_T)^{1/2}x. 
\eeas
When $\frac{2F^2x}{\constr (\log\ttT)r_T^k\ttT}>1$, 
\beas
\frac{x^2}{2\constr (\log\ttT)F^{-2}\tth^{-1}\ttT}
\psi\bigg(\frac{2F^2x}{\constr (\log\ttT)r_T^k\ttT}\bigg)
&=&
\frac{xr_T^k}{8\tth^{-1}}
\times
\bigg(2\cdot
\frac{2x}{\constad (\log\ttT)r_T^k\ttT}\bigg)
\psi\bigg(\frac{2x}{\constad (\log\ttT)r_T^k\ttT}\bigg)
\nn\\&\geq^{(\ref{202412132308})}&
\frac{\constu r_T^kx}{8\tth^{-1}}
\nn\\&\geq^{(\ref{0704012226})}&
\frac{\constu}{8\tth^{-1}}\ttT^{-1/2}(\log\caln_T)^{1/2}x.
\eeas
So we obtained (\ref{202412141033}).

We have
\bea\label{0610061322}&&
\int_{\constn \ttT^{1/2}(\log \ttT)(\log\caln_T)^{1/2}}^\infty {\ared2x}P\big[\big|\bbM_T^k\big|\geq x,{\ared\>\Omega(k,T)}\big]dx
\nn\\&\leq^{(\ref{0704112009})}&
\int_{\constn \ttT^{1/2}(\log \ttT)(\log\caln_T)^{1/2}}^\infty 
4x
\exp\bigg[-\frac{x^2}{2\constr F^{-2}\tth^{-1}\ttT\log\ttT}
\psi\bigg(\frac{2F^2x}{\constr r_T^k\ttT\log\ttT}\bigg)\bigg]dx
\nn\\&\leq^{(\ref{202412141033})}&
\int_{\constn \ttT^{1/2}(\log \ttT)(\log\caln_T)^{1/2}}^\infty 4x
e^{-\consts \ttT^{-1/2}(\log\caln_T)^{1/2}x}dx. 
\eea
}
Applying Lemma \ref{0606231547} in the case where {\bred$q=1$, $p=1$ and $C=\consts \ttT^{-1/2}(\log\caln_T)^{1/2}$}, we obtain the estimate 
{\bred 
\bea\label{0610061323} &&
\int_{\constn \ttT^{1/2}(\log \ttT)(\log\caln_T)^{1/2}}^\infty 4x
e^{-\consts \ttT^{-1/2}(\log\caln_T)^{1/2}x}dx
\nn\\&\simleq&
\ttT(\log\caln_T)^{-1}
\exp\bigg(-\half\consts\constn (\log \ttT)(\log\caln_T)\bigg)
\nn\\&\simleq&
\ttT\caln_T^{-1}
\eea
}
uniformly in $k$, 
{\bred if we take a sufficiently large $\constn$. }%
{\bred 
\begin{lemma}\label{0704081615}
Suppose that 
{\bred
$\ttT\geq3\vee\log\caln_T$ and $\caln_T\geq2$. 
}%
Then, there exists a constant $\constt$ such that 
\bea\label{0610061324}
E\big[\big|\bbM_T\big|^2\big]
&\leq&
\big[\constt \ttT^{1/2}(\log \ttT)(\log\caln_T)^{1/2}\big]^2. 
\eea
The constant $\constt$ is 
depending on $\gamma$, $\tth$, $\big\||\lambda^*|\big\|_\infty$, $\sfd_N$ and $\constzz$, but 
not on $T\in\bbT$. 
\end{lemma}
\proof
From (\ref{0610061321}), (\ref{0704042130}) and (\ref{0610061323}), it is concluded that 
(\ref{0610061324}) holds 
for some constant $\constt$. 
\qed\halflineskip
}
}
\begin{en-text}
\beas 
E\big[\big|\bbM_T^{\tt k}\big|^2\big]
&\leq&
\int_0^\infty 2xP\big[\big|\bbM_T^{\tt k}\big|\geq x\big]dx
\nn\\&\leq&
2\caln_T\sup_k\int_0^\infty xP\big[\big|\bbM_T^k\big|\geq x\big]dx
\eeas

When $T$ is large, 
\beas 
\int_0^\infty xP\big[\big|\bbM_T^k\big|\geq x\big]dx
&\leq&
2\int_0^\infty x
\exp\bigg[-\frac{x^2}{2\constr T}
\psi\bigg(\frac{2F\sfd_N\ol{r}_T^kx}{\constr(\ol{r}_T^k)^2T}\bigg)\bigg]dx
\nn\\&\leq&
(\constn)^2T(\log T)^2\log\caln_T
\nn\\&&
+2\int_{\constn T^{1/2}(\log T)(\log\caln_T)^{1/2}}^\infty
xe^{-\consts\ol{r}_T^kx}dx
\eeas
since 
\beas 
\frac{x^2}{2\constr T}\psi\bigg(\frac{2F\sfd_N\ol{r}_T^kx}{\constr(\ol{r}_T^k)^2T}\bigg)
&\geq&
\frac{x^2}{2\constr T}\frac{\constr\ol{r}_T^kT}{2F\sfd_Nx}
\log\frac{2F\sfd_Nx}{\constr\ol{r}_T^kT}
\nn\\&\geq&
\consts\ol{r}_T^kx
\eeas
for all $x\geq\constn (\log T)T^{1/2}(\log\caln_T)^{1/2}$, where $\consts$ is a constant depending on $F$, $\sfd_N$, 
$\constr$ and $\constn$. 
\end{en-text}

\subsection{Proof of Theorem \ref{0607170208}}

Now we go back to estimation of $R_T$. 
We will basically follow the line of the proof of Theorem 1 in Schmidt-Hieber \cite{schmidt2020nonparametric} 
for the i.i.d. case in the nonparametric regression. 

{\bred 
By choosing a large constant $\constq$, 
it suffices to show that the inequality (\ref{202412220016}) holds 
for sufficiently large $\ttT$, since $\log\caln_T>0$ by the assumption $\caln_T\geq2$, and 
$R_T/(1+F^2)$ is bounded. 
On the other hand, the condition $\ttT\geq \xi(\log\ttT)^2\log\caln_T$ verifies $\ttT\geq3\vee\log\caln_T$ for large $\ttT$. 
Therefore, the assumptions in Lemmas \ref{0703241528} and \ref{0704081615} are satisfied when $\ttT$ is large. 
}

From (\ref{0607160122}), (\ref{0607160130}) and {\ared(\ref{0703240439})}, we obtain 
{\bred 
\bea\label{0607160125}
R_T-R^e_T
&\leq&
\big|R_T-R^e_T\big|
\nn\\&\leq&
\big|\Phi^{(\ref{0512231520})}_T\big|+\big|\Phi^{(\ref{0512231521})}_T\big|
\nn\\&\leq&
{\bred\constzz^2}
F\ttT^{-1}R_T^{1/2}\big\{\big(E[|\bbL_T|^2]\big)^{1/2}+\big(E[|\bbM_T|^2]\big)^{1/2}\big\}
\nn\\&&
+F\ttT^{-1}\bigg\{(\ttT^{-1}(\log \ttT)^{\bblue2}\log\caln_T)^{1/2}+\delta\bigg\}\big(E[|\bbL_T|]+E[|\bbM_T|])
+{\ared2}{\colorr\delta}. 
\eea
}
%
Solving the quadratic inequality (\ref{0607160125}) in $x=R_T^{1/2}$ to obtain
\footnote{We obtain an estimate taking the form $x\leq2^{-1}(A+\sqrt{A^2+4B})\leq A+\sqrt{B}$. 
It gives $x^2\leq2A^2+2B$. 
}
{\bred 
\bea\label{0607160152}
R_T 
&\leq&
2\big[{\bred\constzz^2}F\ttT^{-1}\big\{\big(E[|\bbL_T|^2]\big)^{1/2}+\big(E[|\bbM_T|^2]\big)^{1/2}\big\}\big]^2
\nn\\&&
+2\bigg[R_T^e+F\ttT^{-1}\bigg\{(\ttT^{-1}(\log \ttT)^{\bblue2}\log\caln_T)^{1/2}+\delta\bigg\}(E[|\bbL_T|]+E[|\bbM_T|])+2\delta\bigg],
\nn\\&&
\eea
and hence, for some constant $\constp$ {\ared depending on $\tth$ (for the representation in $T$)} {\bred and $\constzz$}, 
\bea\label{0607160202}
R_T 
&\leq&
2R_T^e+\constp(1+F^2)\bigg[\ttT^{-1}(\log \ttT)^2\log\caln_T+\delta\bigg]
\eea
}
{\bred if $\ttT$ is sufficiently large and $\ttT\geq3\vee\{\xi(\log\ttT)^2\log\caln_T\}$ and $\caln_T\geq2$, }%
from Lemmas \ref{0703241528} and \ref{0704081615}.  
{\bred The condition $\ttT\geq \xi (\log \ttT)^2\log\caln_T$ is used for showing the boundedness of 
$\ttT^{-1}E[|\bbL_T|]$ and $\ttT^{-1}E[|\bbM_T|]$. }


By (\ref{0607160258}), 
\beas
R_T^e
&=& 
E\big[\calr_T]
\nn\\&=&
T^{-1}
E\bigg[-\int_0^T\big\{\wh{a}_T(X_t)-a^*(X_t)\big\}\cdot dN_t+\int_0^T\big\{\wh{b}_T(X_t)-b^*(X_t)\big\}dt\bigg]
\nn\\&=&
T^{-1}E\big[\Psi_T(\wh{a}_T,\wh{b}_T)-\Psi_T(a^*,b^*)\big]
\nn\\&=&
T^{-1}E\big[\Psi_T(\wh{a}_T,\wh{b}_T)-\Psi_T(a,b)\big]
+T^{-1}E\big[\Psi_T(a,b)-\Psi_T(a^*,b^*)\big]
\nn\\&\leq&
\Delta_T+{\tt h}^{-1}E\big[\Psi_{\tt h}(a,b)-\Psi_{\tt h}(a^*,b^*)\big]
\eeas
for any $(a,b)\in{\mathfrak F}_T$.
Therefore, 
\bea\label{0607160300}
R_T^e
&\leq& 
\Delta_T+\inf_{(a,b)\in{\mathfrak F}_T}{\tt h}^{-1}E\big[\Psi_{\tt h}(a,b)-\Psi_{\tt h}(a^*,b^*)\big]
\eea

From (\ref{0607160202}) and (\ref{0607160300}), we obtain Theorem \ref{0607170208}.

\begin{lemma}\label{0606231547}
For positive numbers $p,q,C$ and $B$, let 
\beas 
I(p,q,C,B) &=& \int_B^\infty y^q  e^{-Cy^p}dy.
\eeas
Suppose that 
\bea\label{0606231553}
\frac{q+1}{p}\leq k
\eea
for some number $k$. 
Then there exists a constant $c_k$ depending only on $k$ such that 
\bea\label{0606231555}
I(p,q,C,B) 
&\leq&
c_k\>p^{-1}C^{-\frac{1+q}{p}}\exp\bigg(-\half CB^p\bigg)
\quad\text{whenever}\quad CB^p\geq1. 
\eea
\end{lemma}
\proof
We give a proof here for the sake of self-containedness. 
By change of variables, 
\beas 
I(p,q,C,B)
&=&
\int_B^\infty y^q  e^{-Cy^p}dy
\nn\\&=&
p^{-1}C^{-\frac{1+q}{p}}\int_{CB^p}^\infty u^{p^{-1}(q+1)-1}  e^{-u}du. 
\eeas
Suppose that $CB^p\geq1$. 
Since {\bred$p^{-1}(q+1)-1\leq k-1$} 
by (\ref{0606231553}), we have 
\beas 
I(p,q,C,B)
&\leq&
p^{-1}C^{-\frac{1+q}{p}}\int_{CB^p}^\infty u^{k-1}  e^{-u}du. 
\eeas
There exists a constant $C(k)$ such that $u^{k-1}e^{-u}\leq C(k)e^{-u/2}$ for all $u\geq1$. 
Then 
\beas 
I(p,q,C,B)
&\leq&
C(k)p^{-1}C^{-\frac{1+q}{p}}\int_{CB^p}^\infty  e^{-u/2}du
\nn\\&=&
2C(k)p^{-1}C^{-\frac{1+q}{p}}e^{-CB^p/2}
\eeas
This competes the proof. 
\qed\halflineskip

\section*{Acknowledgement}
The authors thank Professor Taiji Suzuki for the valuable discussion. 

\bibliographystyle{spmpsci}      
\bibliography{bibtex-20180615-20191212-20200312-20200418-20201101++++}   

\appendix

\section*{Appendix}

\begin{figure}[h]
	\includegraphics[width=\textwidth]{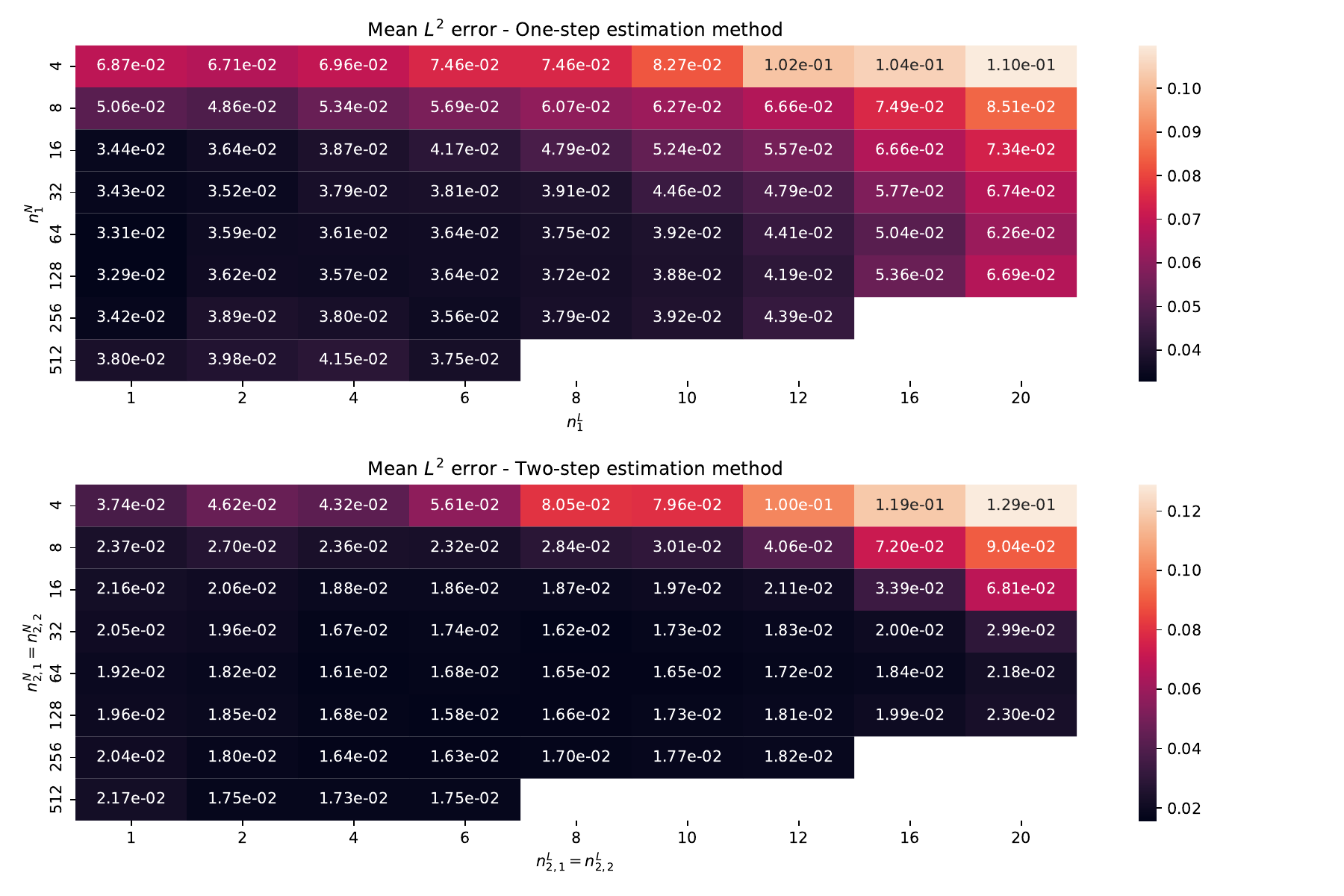}
	\caption{Simulation study --- Heatmap of $L^2$-errors w.r.t the parameters $n^L$ and $n^N$ for both estimation methods.}
	\label{fig:Tuning-L2}
\end{figure}
\begin{figure}
	\includegraphics[width=\textwidth]{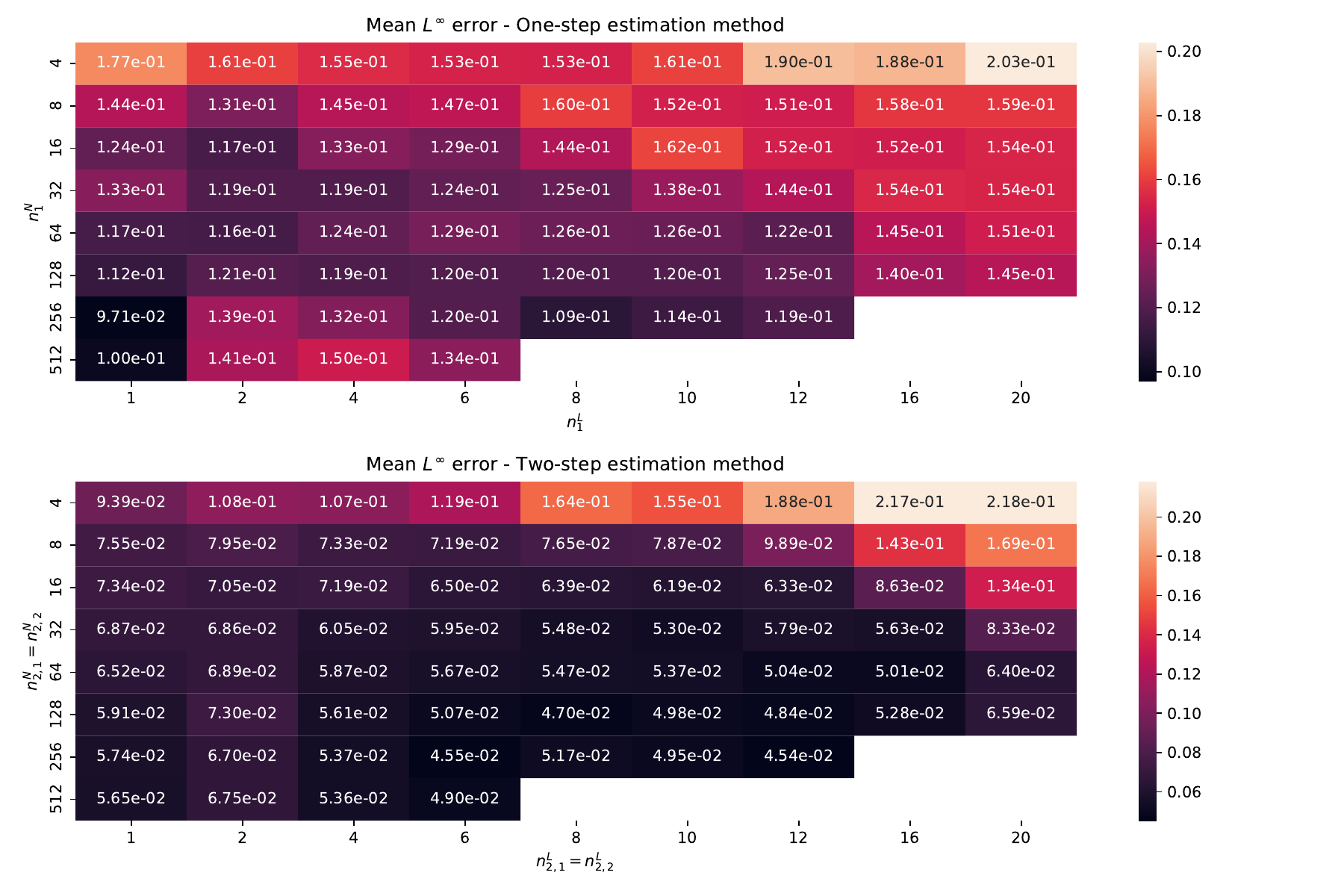}
	\caption{Simulation study --- Heatmap of $L^\infty$-errors w.r.t the parameters $n^L$ and $n^N$ for both estimation methods.}
	\label{fig:Tuning-Linf}
\end{figure}
\begin{figure}
	\includegraphics[width=\textwidth]{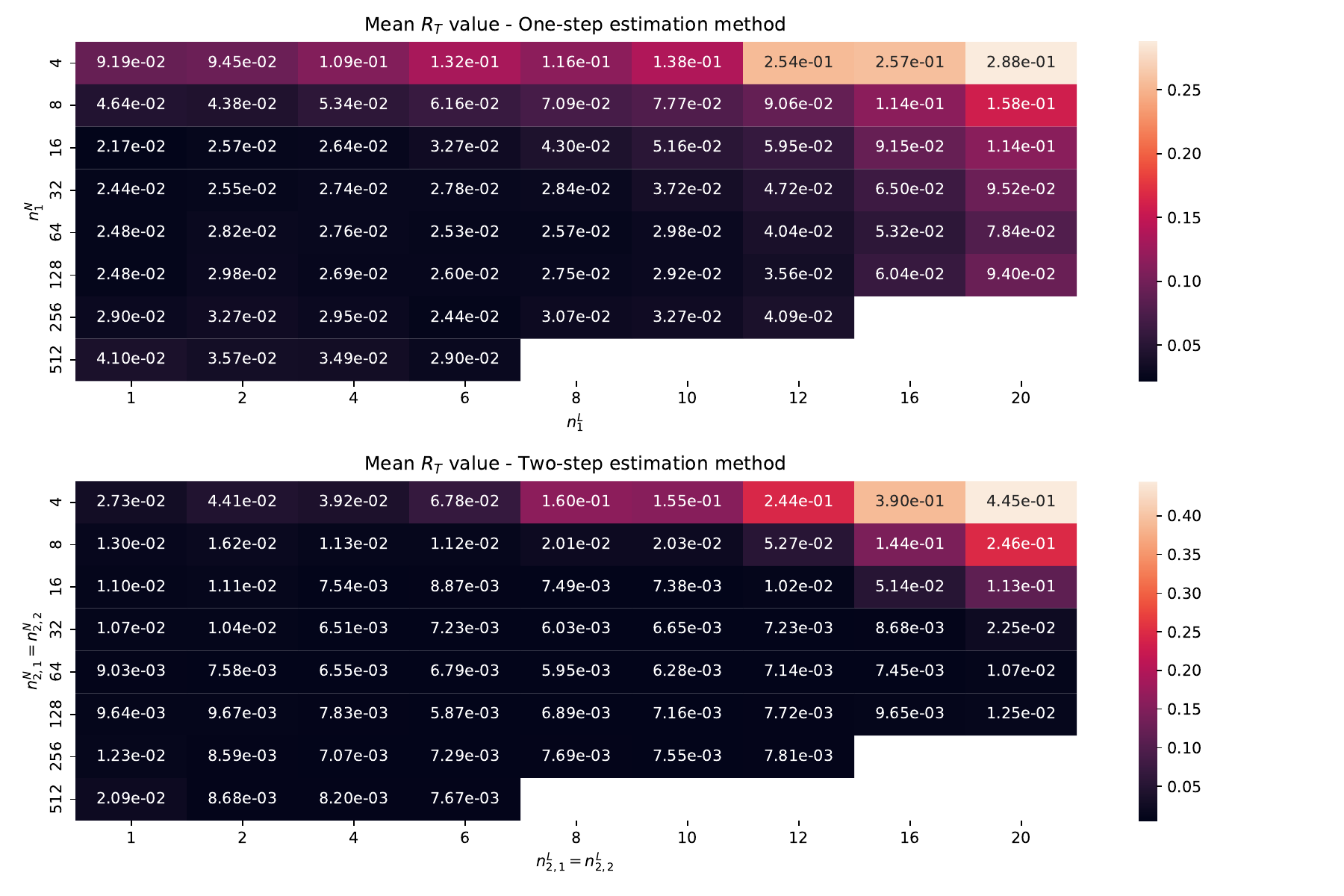}
	\caption{Simulation study --- Heatmap of the empirical values of $\mathcal R_T$ w.r.t the parameters $n^L$ and $n^N$ for both estimation methods.}
	\label{fig:Tuning-RT}
\end{figure}

\end{document}